\documentclass[12pt]{article}
\usepackage{amsfonts, amsthm}
\usepackage{amssymb, amsmath, enumerate, geometry, graphicx, setspace}
\usepackage[T1]{fontenc}
\newtheorem{tw}{Theorem}

\newtheorem{rem}{Remark}

\usepackage[english,polish]{babel}

\begin{document}






\textbf{\Large Dynamical models of dyadic interactions with delay}\\

\begin{center}
\large  Natalia Bielczyk,  Urszula Fory\'s , Tadeusz P\l{}atkowski\\
\small Faculty of Mathematics, Informatics \& Mechanics, Inst. of Appl. Math. \& Mech.,\\
University of Warsaw, ul. Banacha 2, 02-097 Warsaw\\
\end{center}

\begin{abstract}
When interpersonal interactions between individuals are described by the (discrete or continuous) 
dynamical systems, the interactions are usually assumed to be instantaneous: the rates of change of the 
actual states of the actors at given instant of time are assumed to depend on their states at the same time.  
In reality the natural time delay should be included in the corresponding models. 

We investigate a general class of linear models of dyadic interactions with a constant discrete time delay. We prove that in 
such models the changes of stability of the stationary points from instability to stability or vice versa occur for various intervals of the 
parameters which determine the intensity of interactions. The conditions guaranteeing arbitrary number (zero, one ore 
more) of switches are formulated and the relevant theorems are proved. A systematic analysis of all generic cases is carried out. 
It is obvious that the dynamics of interactions depend both on the strength of reactions of partners on their own states as well as on the partner's state. Results presented in this paper suggest that the joint strength of the reactions of partners to the partner's state, reflected by the product of the strength of reactions of both partners, has greater impact on the dynamics of relationships than the joint strength of reactions to their own states. The dynamics is typically much simpler when the joint strength of reactions to the partner's state is stronger than for their own states. Moreover, we have found that multiple stability switches are possible only in the case of such relationships in which one of the partners reacts with delay on their own state. 
Some generalizations to triadic  interactions are also presented.  
\end{abstract}


\newpage

\textbf{\Large Dynamical models of dyadic interactions with delay}\\
\bigskip{}
\section{Introduction}
Interpersonal relationships such as dyadic interactions are very complicated and difficult to predict. One of the possibilities 
of making some predictions is to use mathematical modeling. Mathematical modeling of interactions between humans is a relatively 
new branch of applied mathematics. Typically, in mathematical sociology, reactions of individuals to stimuli are assumed to be 
instantaneous, both in discrete and continuous models, cf. for example~\cite{Bar, Bud, Now, Fel, Got, Lat, Lie, Mee, Rin, Rus, 
Str1, Str2}  and the references therein. In the basic models of relationships between two partners the linear approach can be 
used, cf.~\cite{Fel, Str1, Str2}. It is obvious that such simple models cannot describe the relationships 
precisely. However, they can capture some key features of them. On the other hand, it is well known that all 
natural phenomena are delayed with respect to its stimuli. In the case we study, taking into account the delay 
of reactions can better describe the real dynamics of interpersonal relations because of the natural tendency of some people to analyze their own and/or partner's reactions and act after that, which necessarily takes some time. 

In the simplest approach, the stimuli are the actual states of both partners. Instantaneous reactions correspond formally to the assumption that the speed, the rates of change of the states and the actual states of the partners are considered at the same instant of time. However, because time delays are present in  human interactions, the reactions to stimuli are naturally retarded. Whereas in biological applications the delay has been introduced many years ago, cf.~the Hutchinson model~\cite{hut} which has been proposed in~1948 as a delayed version of the logistic equation, and the models with delays have been already considered in numerous biomathematical literature, cf. for example  \cite{For1, Gop, Hal, murray} and the references therein, the delay in systems modeling human behavior have not received as much attention as it possibly should.  

In this paper we study possible scenarios of the evolution of the linear models in which the constant in time 
discrete delay is introduced. We carry out a systematic study of the models with delay present in different stimuli 
terms corresponding to the states of the partners. Although in general nonlinear models  are more appropriate in 
the description of complicated systems, the mathematical approach we use in the paper allows us to simplify the 
description to the linear one. More precisely, we study small deviations of the partner states from some states 
(which are the equilibrium or stationary states of the system) and check whether  these deviations tend to zero 
(stable case) or not (unstable case). Because the considered deviations are small, we are able to approximate the 
nonlinear dynamics using linear models, as usually in the stability theory for differential equations without delays 
(ODEs) and with delays (DDEs), cf. e.g. ~\cite{Hal, Hal1}.   

In general, different scenarios of influence of a discrete delay on the stationary state of the model are possible. 
\begin{enumerate}
\item The delay can have no effect on the stability of the stationary state. 
\item Sufficiently large delay can destabilize the stationary state which is stable without the delay.
\item The appropriate values of delay can stabilize the stationary state which is unstable  without the delay. 
\end{enumerate}
The second case is the most typical. In this paper we show however, that the third case occurs in the systems which are modeled by two linear DDEs for different intervals of the coefficients  determining the model. 

\subsection{Basic Romeo and Juliet model}
As an introductory example, let us consider the model of the mutual interactions of two persons, Romeo (R) and Juliet (J), without the delay, cf.~\cite{Str1, Str2}.   
Let the functions  $R(t)$ and $J(t)$ denote the state of R and J, respectively. The state of the partner can be interpreted as the satisfaction, the level of happiness, etc. 
from their mutual relation. Positive values reflect the positive attitude (good emotions) towards the partner, or positive experience from the relation, negative values reflect the negative, bad experience. 

In~\cite{Str1,Str2}, the following model has been considered (we cite the author): ''R loves J, J is a fickle lover. The more R loves her, the more J wants to run away and hide. When R gets discouraged and backs off, J begins to find him strangely attractive. R on the other hand tends to echo her: he warms up when she loves him, and grows cold when she hates him.''. In~\cite{Str1,Str2},  this model is described by the system of two linear ODEs 
\begin{equation}
\left\{
\begin{array}{lcr}
\dot{R}(t) & =  & aJ(t) ,\\
\dot{J}(t) & =  & - bR(t),
\end{array}
\right. 
\label{eq:2}
\end{equation} 
where  $a$, $b  > 0$, according to the interpretation above (other nontrivial cases lead to divergent solutions). The solutions to Eqs.~\eqref{eq:2} oscillate around the equilibrium at $(0,0)$ which is a center in the phase space $(R,J)$.

\subsection{General R and J model}
Now, consider more general situation. Let R and J have their fixed goals or preferences, ideals of their states in the relation, which we denote by $R^*$ and $J^*$, respectively. Assume that the rate of change of their actual states  depends linearly on two factors: the differences between their own actual states and their own ideals, and the difference between the actual states of the  partners. Such model has been studied in~\cite{Fel}, and reads
\begin{equation}
\left\{
\begin{array}{lclcl}
\dot{R}(t) & = & c[R^* - R(t)] & + & d[J(t) - R(t)] ,\\
\dot{J}(t) & = & e[J^* - J(t)] & + & f[R(t) - J(t)] .\\
\end{array}
\right. 
\label{eq:248}
\end{equation} 
The model enables to assign clearly personal features of the persons engaged and compute the evolution of the system easily. The coefficients $c$, $d$, $e$, $f$ can be positive or negative, and describe different psychological traits of the partners, for example contrarian, cooperator etc. We refer the reader to~\cite{Fel} for the discussion of different aspects and interpretation of the model. 

We are interested in the temporary evolution of the states of the partners. 
Let us denote by $(R_{eq}, J_{eq})$ the stationary state of this linear inhomogeneous system. The stability properties of this state do not change under the linear transformation of the independent variables. Substituting  
$$ r(t):=R(t)-R_{eq}, \ \ \quad  j(t):=J(t)-J_{eq} ,$$ 
we obtain the homogeneous system 
\begin{equation}
\left\{
\begin{array}{lcl}
\dot{r}(t) & = & a_{11}r(t) + a_{12}j(t) ,\\
\dot{j}(t) & = & a_{21}r(t) + a_{22}j(t) ,
\end{array}
\right. 
\label{eq:8}
\end{equation} 
with the equilibrium $(0,0)$, and general coefficients $a_{kl} \in \mathbb{R}$, $k, l\in \{1, 2\}$.   

The conditions for stability of the equilibrium for this system can be found in textbooks on ODEs, cf.~\cite{Hal1}. In terms of the dyadic interactions they can be interpreted as conditions of stable relations. They strongly depend on the signs of the coefficients $a_{kl}$ and the relations between them. 
The choice of signs of the coefficients $a_{kl}$ specifies different psychological traits of the partners. The case $a_{11} >0$, $a_{12} >0$ can be referred to as ''eager beaver'', cf.~\cite{Str2}: ''R gets excited by J love for him, and is further spurred on by his own affectionate feelings for her.'' The case $a_{11}<0$, $a_{12}>0$ is referred to as ''cautious lover''.  The case $a_{11}>0$, $a_{12}<0$ could be called ''narcissus'', and $a_{11}<0$, $a_{12}<0$ -- the doubly cautious, or ''stoic lover''.  It should be noticed that the model is symmetric with respect to its variables, therefore the same interpretation as we have given for R is also valid for J. 

\subsection{Delayed R and J model}
The models of dyadic interactions described by Eqs.~\eqref{eq:8} do not take into account possible delays in the reactions of the partners to their instantaneous states.  In other words, the rate of change of the personal state  depends only on the actual state of this person and that of the partner. 

In this paper we consider situations in which the delays of the reactions of the partners on their own states 
and/or the states of their partners are taken into account. In general, the delay can have various effects on the 
dynamics of the partnership, for example we show that delay can help two cautious lovers to 
stay in a stable relationship although their relationship would lose stability without it. 

We consider the case when the contributions to the rates of changes of the states $r(t)$ and/or $j(t)$ can be influenced by the past states of one or two of the partners, that is
\begin{equation}
\left\{
\begin{array}{lcl}
\dot{r}(t) & = & a_{11}r(t-\tau_{11}) + a_{12}j(t-\tau_{12}) ,\\
\dot{j}(t) & = & a_{21}r(t-\tau_{21}) + a_{22}j(t-\tau_{22}) ,
\end{array}
\right. 
\label{eq:ogolne}
\end{equation} 
where $\tau_{kl} \geq 0$, $k, l\in \{1, 2\}$, are the delays of reactions, that is for $\tau_{kk}>0$, $k\in \{1, 2\}$ the person reacts with delay on his/her own state while for $\tau_{kl}>0$, $k\ne l$, $k, l \in \{1, 2\}$, the delay is present in the reaction to the partner state. It is obvious that if $\tau_{kl}=0$ for some $k, l\in \{1, 2\}$, then this type of reaction is instantaneous. 

In general, it is difficult to study the dynamics of Eqs.~\eqref{eq:ogolne} for different values of delays $\tau_{kl} \ne \tau_{mn}$ for $kl\ne mn$. However, the case when the contribution or contributions from the past states is delayed by the same time interval can be analyzed almost completely. Thus, we consider Eqs.~\eqref{eq:ogolne} with one delay $\tau$ present in some of the reactions terms, while other reactions 
remain instantaneous. 
Moreover, it turns out that even in such drastically simplified model interesting effects occur, in particular the appropriate delay can stabilize the equilibrium outcome of the evolution of the states of the couple. 
Such effect has been also found in~\cite{Bod} for a model with an additional linear term.

It is well known, cf.~\cite{Coo}, that if the delay is large enough, and the equilibrium is unstable for at least one value of the delay $\tau \ge 0$, then the trajectories of the considered system diverge for delays large enough. In real applications the behavior of the partners for bounded delays is more relevant. Let us consider the arbitrary maximal value of delay $\tau_{\max}$. We can put the following questions.
\begin{enumerate}
\item The equilibrium is unstable if there are no delays. Can a larger delay $0<\tau<\tau_{\max}$ stabilize it? 
\item The equilibrium is stable if there are no delays. If the delay is large enough, $0<\tau<\tau_{\max}$, it has the stabilizing effect. What conditions on the reactive parameters $a_{kl}$ enlarge the interval of stability? 
\item How many stability switches are possible for the situations 1. and 2.?   
\item For given parameters of one partner, what can the second person do in order to stabilize the relationship? When is there a possibility 
to gain stability by changing one's attitude towards themselves and the other person's state, or the relevant delays? 
\end{enumerate}
The last question can also be formulated as: can a relationship be in command of one flexible person, when the partner is stiff in attitudes and level of reflexivity?

\vskip 0.1cm
In the following we address these issues, prove the relevant results and discuss numerical examples.

\subsection{Notions of stability and instability}
Now, we introduce the notion of stability, instability and stability switches.
\begin{itemize}
\item We say that the equilibrium is stable, when solutions to Eqs.~\eqref{eq:ogolne} tend to this state in time. 
\item We say that the solution is unstable, when solutions to Eqs.~\eqref{eq:ogolne} diverge. 
\item We say that the system  is stable when its unique equilibrium is stable.
\item We say that the system is unstable, when the unique equilibrium is unstable.
\item We say that the stability switch occurs if there is the change of stability of equilibrium from unstable to stable or vice versa, 
caused by a change of a parameter of the system. In this case typically a Hopf bifurcation occurs.  
\end{itemize}
For precise definitions we refer to any textbook in ODEs, cf.~\cite{Hal1}, or DDEs, cf.~\cite{Hal}.

It should be noticed that for equations of the form~\eqref{eq:ogolne}  there is a unique equilibrium $(0,0)$, while 
for the more complicated, especially non-linear models, there can be more equilibriums. In general, we can also 
consider another type of stability, the Lapunov stability, when solutions starting near equilibrium stay near it 
but not necessarily tend to it. Typically, Lapunov stability is non-generic, which means that the small change of the model parameters leads to stability or instability of equilibrium in the sense described above. In the paper we mainly focus on the generic cases, as typical in 
analysis of the models of natural phenomena. This approach reflects the fact that in reality the models parameters do not 
stay exactly the same all the time but can change in time and therefore, non-generic cases are very difficult to observe. 

The terms on the right-hand side (rhs) of Eqs.~(\ref{eq:ogolne}) are addressed as {\it{stimuli terms}} or reaction terms. 

\vskip 0.1cm
In the next three sections we consider the models with the delay placed in one, two, three and four terms. In section V we introduce a model of triadic interactions. In section VI we conclude and discuss some open problems.  
In Appendix the proofs of theorems stated in the paper and some other notes on the topic are presented.

\section{Single delay}
First we consider the relationships in which only one of the stimuli terms is delayed. There are two families of such models.   
In the first one the evolution of R's state depends on his state in the past, that is instead of the general system~\eqref{eq:ogolne} we study
\begin{equation}
\left\{
\begin{array}{lcccl}
\dot{r}(t) & = & a_{11}r(t-\tau) & +  & a_{12}j(t) ,\\
\dot{j}(t) & = & a_{21}r(t) & + & a_{22}j(t) .
\end{array}
\right. 
\label{eq:12}
\end{equation} 
In other words, R reacts immediately to the state of J, whereas is reflexive to his own states, and reacts towards them only 
after the delay $\tau$. This is a reasonable model when it comes to a person who tends to analyze the 
environment rather than themselves and their own state affects them only after some deliberation. 
J is more self-sensitive in this model, she tends to update her attitude on the basis of temporary 
level of private satisfaction. Mathematical properties of this model are described by Theorems~1 and~2 below. 

In the second model the evolution of the R's state depends on the J's delayed stimulus, 
whereas his reaction to his own state is immediate, occurs without delay. It is described by the following form of Eqs.~\eqref{eq:ogolne}  
\begin{equation}
\left\{
\begin{array}{lcccc}
\dot{r}(t) & = & a_{11}r(t) & + & a_{12}j(t-\tau) ,\\
\dot{j}(t) & = & a_{21}r(t) & + & a_{22}j(t) .
\end{array}
\right. 
\label{eq:23}
\end{equation} 
This system describes the situation when R is impulsive about his internal states and follows his 
satisfaction level consciously by reacting to it immediately, but at the same time he is cautious about J and reacts to 
her levels of satisfaction with delay. This may result from carefulness as well as lack of knowledge about 
her feelings. Mathematical properties of this model are described by Theorem 3 below. 

It is obvious that exactly the same models can be considered when we replace the model variables, that is the roles of R and J are replaced. We omit details for such models because its dynamical properties are the same as those presented below. 
We look for stability switches of the equilibrium for both models~(\ref{eq:12}) and~(\ref{eq:23}). 
We remind, cf. for example~\cite{Coo}, that if at least one stability switch occur, then $(0,0)$ is unstable for delays large enough. 
Note also that Eqs.~(\ref{eq:12}) and~(\ref{eq:23}) are stable for $\tau = 0$ if

\begin{equation} \label{warunek0}
a_{11}+a_{22} < 0 \quad  \text{and} \quad  a_{11}a_{22} - a_{12}a_{21} >0. 
\end{equation} 
Inequalities~\eqref{warunek0} can be also expressed by the properties of the matrix 
$\mathbf{A}=\left( \begin{array}{cc}
a_{11} & a_{12} \\ a_{21} & a_{22}
\end{array}\right)$ which is formed by the model coefficients, namely
\[
\text{tr} \mathbf{A}=a_{11}+a_{22}<0 \quad \text{and} \quad \det \mathbf{A} =a_{11}a_{22}-a_{12}a_{21}>0 .
\] 
Note also that for $\text{tr} \mathbf{A}=0$ there is Lapunov stability if $\det \mathbf{A}>0$. 

In Appendix we prove the following theorems.
\begin{tw} \label{tw1} 
For $|a_{12}a_{21}|<|a_{11}a_{22}|$
\begin{enumerate}
\item  if the system~(\ref{eq:12}) is unstable for $\tau=0$, then it is unstable for all $\tau \ge 0$;
\item if the system~(\ref{eq:12}) is stable for $\tau=0$, then there is an unique stability switch (for which there is the Hopf bifurcation).
\end{enumerate}
\end{tw}
Thus, for $|a_{12}a_{21}|<|a_{11}a_{22}|$ there is at most one stability switch for all admissible reactive parameters $a_{kl}$. 

This result suggests that if the joint strength of reactions of partners to their own states, reflected by the product of both strength $|a_{11}a_{22}|$, is larger than the joint strength of reactions 
to partners' states, there is no 
stability at all or the system loses stability once for all. The conclusion is that, as far as such a relationship was stable and in a result of slight moderation of parameters (which usually accompanies personal development of partners) starts to shiver, the delayed person should just 
accelerate reactions and the stability will come back to the system. 

Note that the behavior of the system depends continuously on the parameters of the model, 
which implies that decrease in delay will cure the relationship even if it loses stability as a result of 
a slight change in attitude of partners (i.e. some of the parameters $a_{kl}$, $k, l \in \{1, 2\}$). This situation corresponds to reality and natural development of people because it is common to moderate attitudes and reactions in long periods of time.  

Similar interpretation can be given to the results for the model~(\ref{eq:23}), cf.~Theorem~\ref{tw3} below.

For $|a_{12}a_{21}|>|a_{11}a_{22}|$ we define  
\begin{equation}
\Delta := (a_{11}^2 + a_{22}^2)^2 + 4a_{12}a_{21}(a_{22}^2 - a_{11}^2). 
\label{eq:18000}
\end{equation}
For $\Delta < 0$, any of the systems considered in this paper has no stability switches, which is proved in Appendix.
The case $\Delta = 0$ is non-generic and therefore we omit it in our analysis.
For $\Delta > 0$ and $i\in \{0, 1\}$ we denote 
\begin{equation}
y_i := \frac{1}{2}\left(a_{11}^2 - a_{22}^2 - 2a_{12}a_{21} \pm \sqrt{\Delta}\right), \quad \omega_i:=\sqrt{y_i} \ \ \text{for} \ y_i > 0,
\label{eq:18001}
\end{equation}
\begin{equation}
a_i:=\frac{a_{12}a_{21}a_{22}}{a_{11}(\omega_i^2 + a_{22}^2)} .
\label{eq:18002}
\end{equation}   

\begin{tw} \label{tw2} 
If $|a_{12}a_{21}| > |a_{11}a_{22}|$, then for all natural $n \in \mathbb{N}$, there exist intervals of ${a_{kl}}$, $k, l\in \{1, 2\}$, such that Eqs.~\eqref{eq:12} has $n$ stability switches.
The set of sufficient conditions reads
\begin{equation}
\Delta > 0, \quad 
\quad a_{11}^2 – a_{22}^2 - 2a_{12}a_{21} > 0
\label{eq:tencozwykle}
\end{equation}
and for $i\in \{1, 2, ... ,2(\lceil \frac{n}{2} \rceil -1)\}$ 
\begin{equation}
\left\{
\begin{array}{l}
\frac{1}{\omega_1} \arccos a_0 + \frac{(2i - 2) \pi}{\omega_1}   <   \frac{1}{\omega_0} \arccos a_1 + \frac{(2i - 2) \pi}{\omega_0} , \\
\frac{1}{\omega_0} \arccos a_1 + \frac{(2i - 2) \pi}{\omega_0}   <   \frac{1}{\omega_1} \arccos a_0 + \frac{(2i) \pi}{\omega_1} .
\end{array}
\right. 
\label{eq:6000}
\end{equation}
\end{tw}
The appearance of a stability switch  implies instability of the system for $\tau \rightarrow  + \infty.$ 
Thus, in the generic case, if the system is unstable for $\tau=0$, the number of switches must be even, while 
if it is stable for $\tau=0$, the number of switches must be odd.
However, we should remember that in reality the delay is bounded, $\tau<\tau_{\max}$, and therefore the number of switches  observed in the interval $[0,\tau_{\max}]$ can be different than those obtained from the analytical investigation when we let $\tau \to +\infty$. 

Hence, for $|a_{12}a_{21}| > |a_{11}a_{22}|$, there may be narrow segments of delay values for which the relationship is stable. This means the partners cannot know what will be the impact of decreasing or increasing the delay on the stability of the relationship, which implies that if the system loses stability, it is no longer easy to regain it on the basis of conscious acting. 
This system describes a relationship in which the joint strength of reactions of partners to their own states is weaker than the joint strength of reactions 
to partners' states. This contributes to a common belief that relationship in which partners are more emotionally focused on partner's states than on their own, are shaky, because it is essential to hold one's own internal source of satisfaction and emotional spine. 

\begin{rem}
We can also consider other types of linear models, e.g. the situation in which the evolution of $r(t)$ is influenced by its current value $r(t)$ and the past value $r(t-\tau)$ (and the partner's actual state $j(t)$)  
\begin{equation}
\left\{
\begin{array}{lccclcl}
\dot{r}(t) & = & a_{11}r(t-\tau) & + & a_{12}j(t) & + & a_{13}r(t) ,\\
\dot{j}(t) & = & a_{21}r(t) & + & a_{22}j(t) && .
\end{array}
\right. 
\label{eq:410}
\end{equation} 
An example of two stability switches for this model has been considered in~\cite{Bod}. 
\end{rem}
We note that for $a_{13} = 0$ the system~(\ref{eq:410}) is identical to~(\ref{eq:12}). Since the stability properties of the system depend continuously on the parameters, there exists an interval around $a_{13}=0$ for which the stability properties of both systems are the same.
Since any number of switches can be obtained for Eqs.~(\ref{eq:12}), the same is true for Eqs.~(\ref{eq:410}) with $a_{13}$ small enough.
We argue it more precisely in Appendix.

\vskip 0.2cm
For the second model~(\ref{eq:23}), in which the rate of change of the R's state depends on the J's 
state in the past, unlike for Eqs.~(\ref{eq:12}), at most one stability switch can occur.  
\begin{tw} \label{tw3}
For $|a_{12}a_{21}|<|a_{11}a_{22}|$ there are no stability switches for the system~(\ref{eq:23}). 

For $|a_{12}a_{21}|>|a_{11}a_{22}|$ 
\begin{enumerate}
\item if the system~(\ref{eq:23}) is unstable for $\tau=0$, then it is unstable for all $\tau \ge 0$;
\item if the system~(\ref{eq:23}) is stable for $\tau=0$, then there is an unique stability switch. 
\end{enumerate}
\end{tw}

\begin{rem}
We note that if we add the delay to the simplest model~(\ref{eq:2}), which is the particular case of~(\ref{eq:12}), 
\begin{equation}
\left\{
\begin{array}{lcc}
\dot{r}(t) & = & a_{12}j(t-\tau) ,\\
\dot{j}(t) & = & - a_{21}r(t), 
\end{array}
\right. 
\label{eq:4}
\end{equation} 
with $a_{12}$, $a_{21}  > 0$, it turns out that, apart from particular initial data, and a measure zero set 
of discrete values of delays, the states of the partners diverge. 
This system can be rewritten as 
$r''(t) = - \omega^2 r(t - \tau)$, where 
$\omega^2 = a_{12}a_{21}$, which only gives stable results for $\tau = 2k\pi/\omega$, $k \in \mathbb{N}$. 
The same is true if both partners react with different delays, $\tau_1$ and $\tau_2$, respectively, because then 
the system can be rewritten as $r''(t) = - \omega^2 r(t-(\tau_1+\tau_2))$, which means it behaves in the 
same manner, but with the delay  $\tau=\tau_1 + \tau_2$.
However, this is stability in the Lapunov sense   and therefore it is not preserved for the parameters near these critical values. 
 Thus, in this case the delay cannot lead to stabilization of the relation. We omit details. 
\end{rem}
Below we give an example of the system with a single delay, for which no  stability switch occurs and the system is 
stable for every value of delay. Let us consider Eqs.~\eqref{eq:12} 
with $\Delta < 0$ and the system is stable for $\tau = 0$, which provides a set of conditions
\begin{equation}
\left\{
\begin{array}{c}
(a_{11}^2 + a_{22}^2)^2 + 4a_{12}a_{21}(a_{22}^2 - a_{11}^2) > 0 ,\\
a_{11} + a_{22} < 0 ,\\
a_{11}a_{22}- a_{12}a_{21} > 0 .
\end{array}
\right. 
\label{eq:8003}
\end{equation}
Conditions~(\ref{eq:8003}) are satisfied e.g.~for the coefficients 
$a_{11} = a_{12} = 1$, $a_{21} = a_{22} = -2$. In this example R is an eager beaver, and J is a 
stoic/double cautious lover, and both of them are much more reactive for partner's state than their own. 
The delay does not affect the stability of the system, it affects however the range of time we have to wait for the 
system to stabilize. This means that in such a relationship R can be reflexive to any extent, and the relationship is still 
stable. 

\vskip 0.2cm
In figures below we present numerical examples of the stability switches and the Hopf Bifurcation (HB) for 
different choices of the parameters of the first model~(\ref{eq:12}). In all figures the coefficients 
$a_{kl}$ of the system~(\ref{eq:12}) are represented as the coordinates of the vector $[a_{11},a_{12},a_{21},a_{22}]$. 
The initial data are  $r(t)\equiv 1$ and $j(t)\equiv 1$, $t\in [-\tau,0].$ 
\begin{figure}[ht]
\centering
\includegraphics[width=5.5cm]{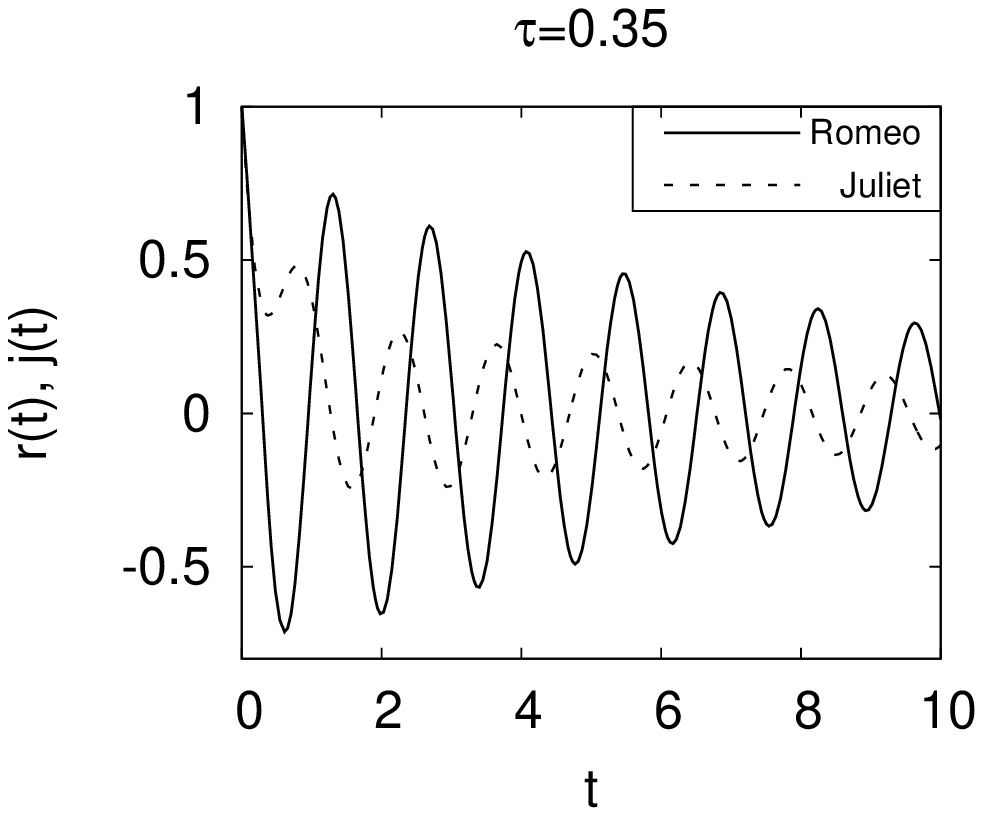}\includegraphics[width=5.5cm]{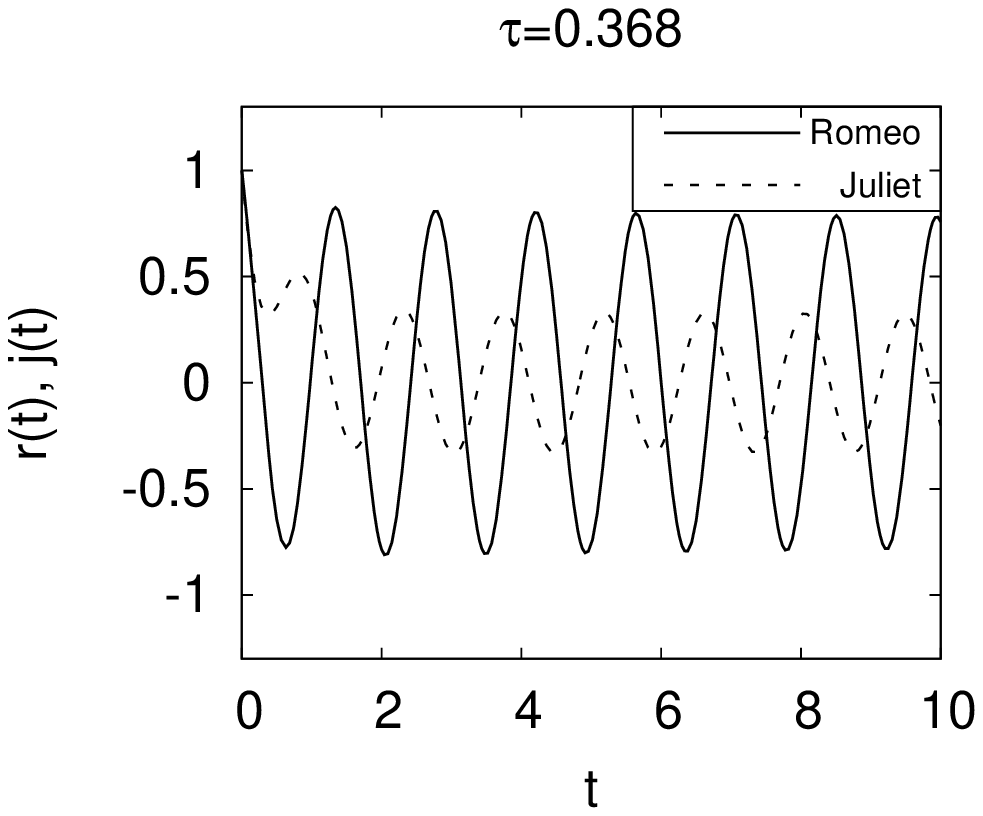}\includegraphics[width=5.5cm]{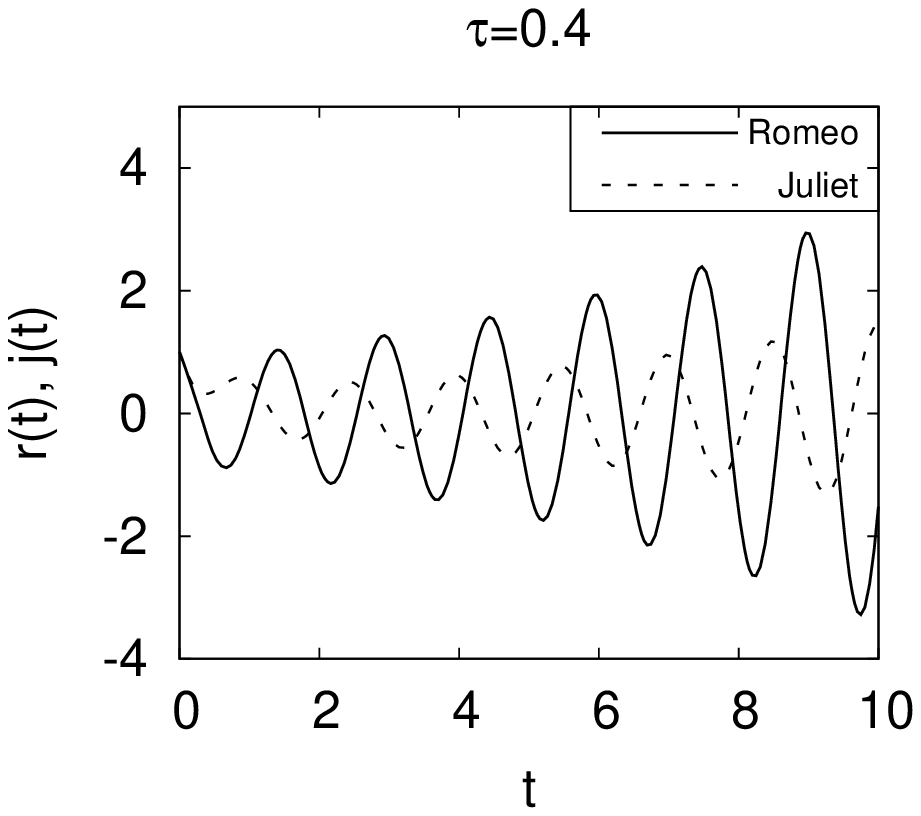}
\caption{Time evolution of solutions $r(t)$, $j(t)$ to Eqs.~(\ref{eq:12}) for $[-4, 1, -2, -2]$ around the unique stability switch 
$\tau \approx 0.3687$. 
Left figure: $\tau=0.35$, stable; middle: $\tau \approx 0.368$, HB; right: $\tau= 0.4$, unstable}
\label{fig0}

\centering
\includegraphics[width=5.5cm]{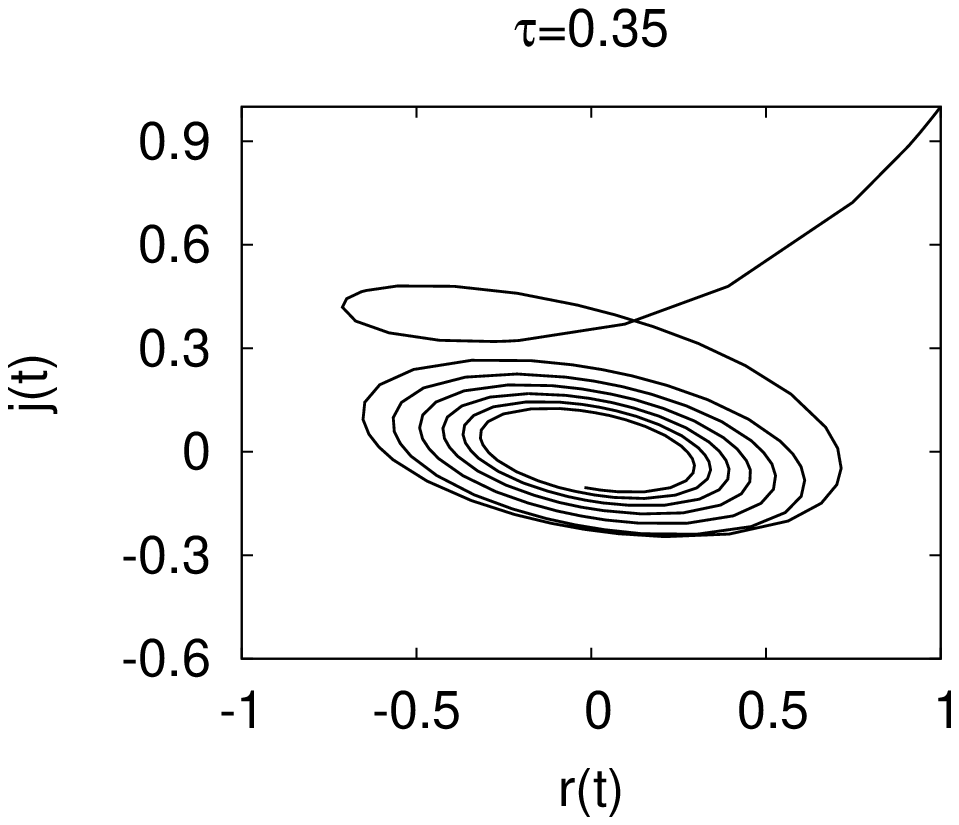}\includegraphics[width=5.5cm]{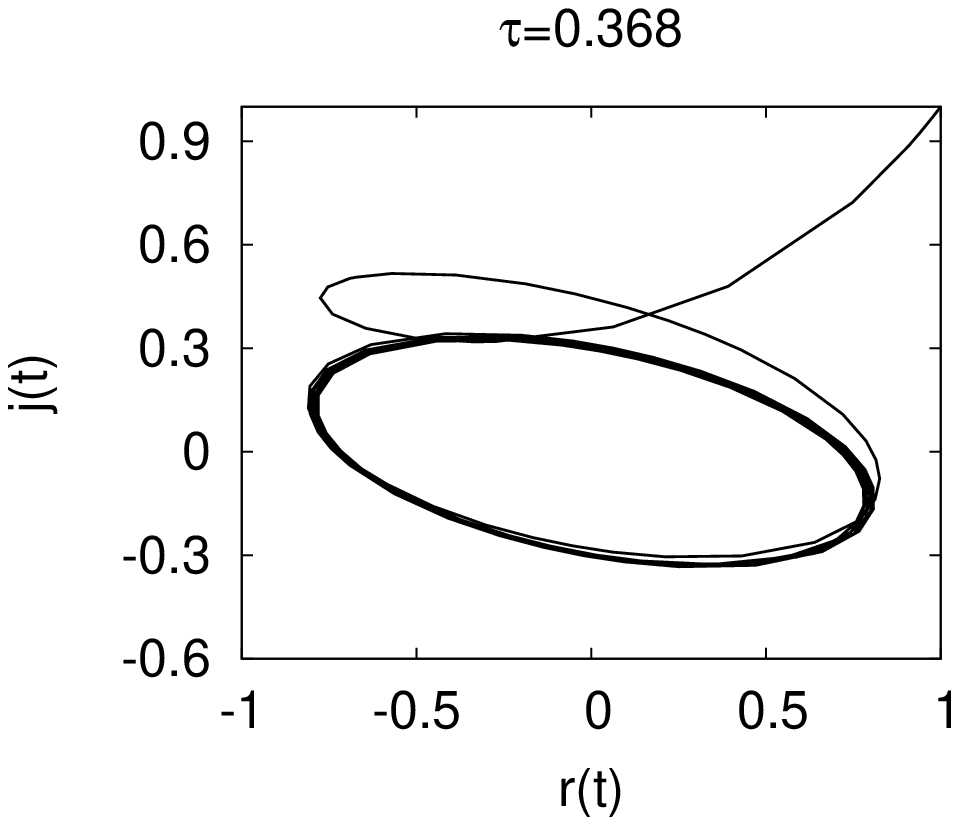}\includegraphics[width=5.5cm]{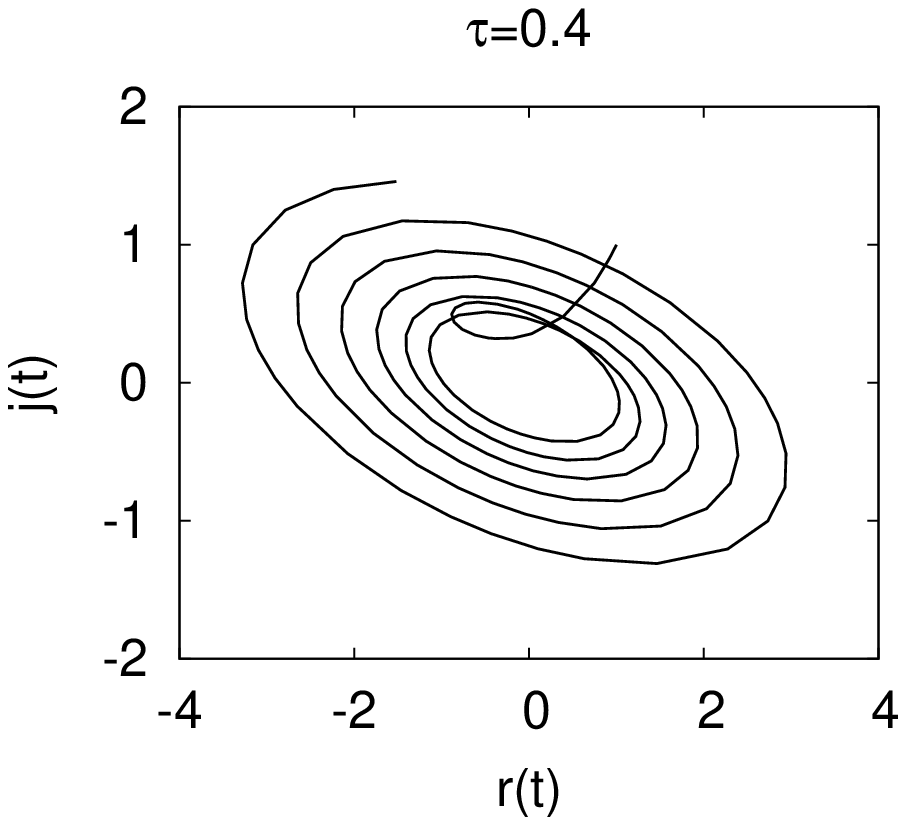}
\caption{Trajectories  of Eqs.~(\ref{eq:12}) in the space $(r, j)$ for $[-4,1,-2,-2]$ around the unique stability switch 
$\tau \approx  0.3687$. Left figure: $\tau=0.35$, stable; middle: $\tau \approx 0.368$, HB; right: $\tau= 0.4$, unstable}
\label{fig1}
\end{figure}

In Fig.~\ref{fig0} we show the trajectories $r(t)$, $j(t)$ for delay values close to the unique stability switch for the reactive 
parameters $a_{kl}$ corresponding to p.~2 of Theorem~\ref{tw1}, that is 
$[a_{11}, a_{12}, a_{21}, a_{22}]=[-4, 1, -2, -2]$. For delays below the switch value, the equilibrium is stable, for larger 
delays is unstable. The left diagram shows the trajectory which converge to the equilibrium. The middle one visualizes 
the limit cycle appearing as the result of Hopf bifurcation, and the right one presents the diverging trajectory. 
In Fig.~\ref{fig1} we show the trajectories in the  space $(r, j)$, corresponding to above 
values of the delay. Similar behavior is obtained e.g. for the reactive parameters $[-4,-2,-5,-6]$, $[-4,4,4,-9]$, $[-5,-2,2,-5]$. As it can be seen at least two of $a_{kl}$, that is $a_{11}$ and $a_{22}$ are negative in this case.

\begin{figure}[ht]
\centering
\includegraphics[width=5.5cm]{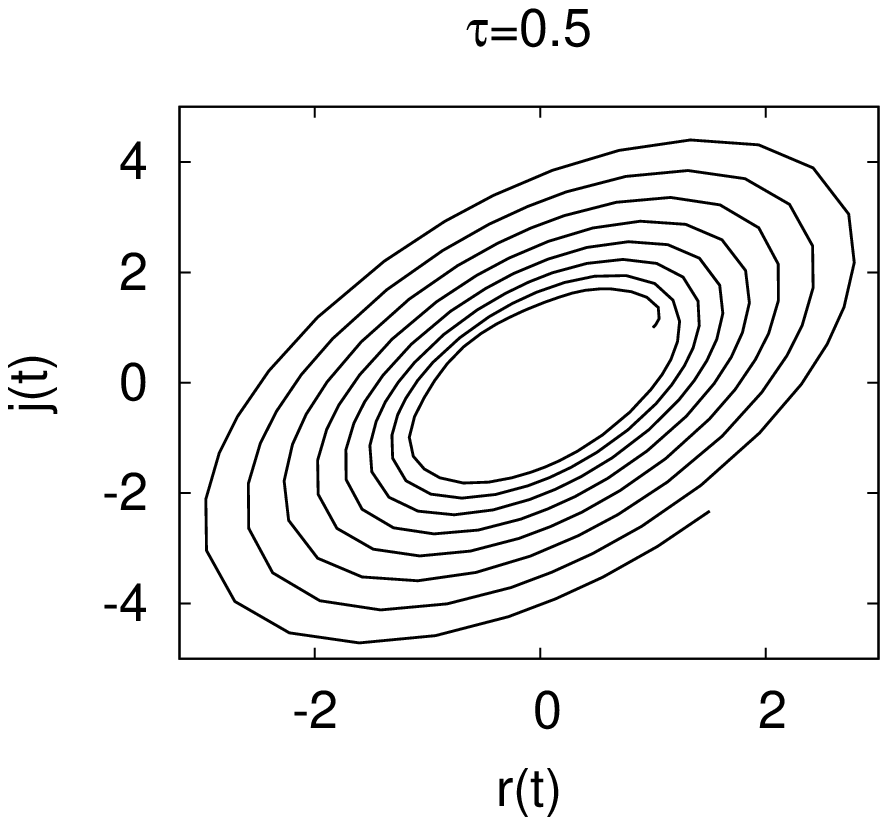}\includegraphics[width=5.5cm]{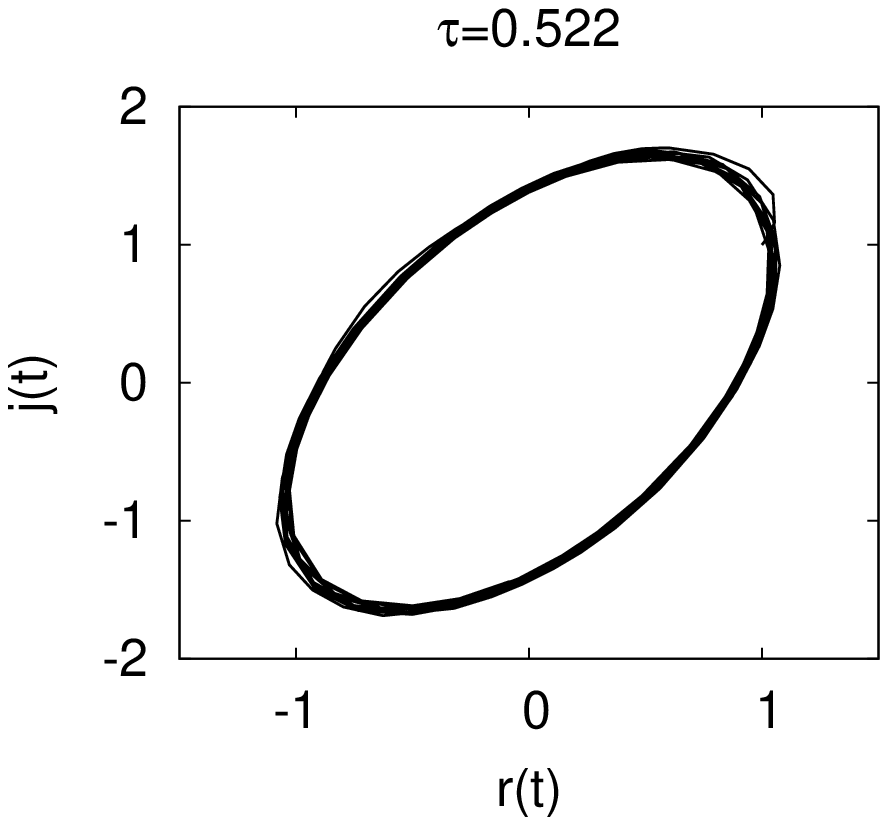}\includegraphics[width=5.5cm]{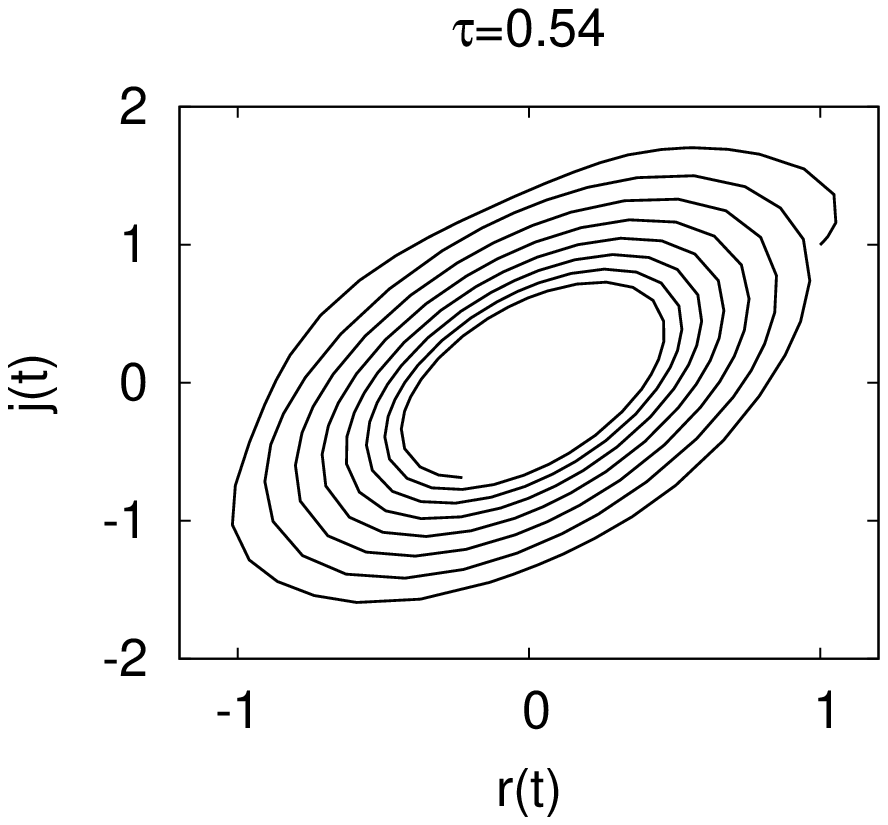}

\centering
\includegraphics[width=5.5cm]{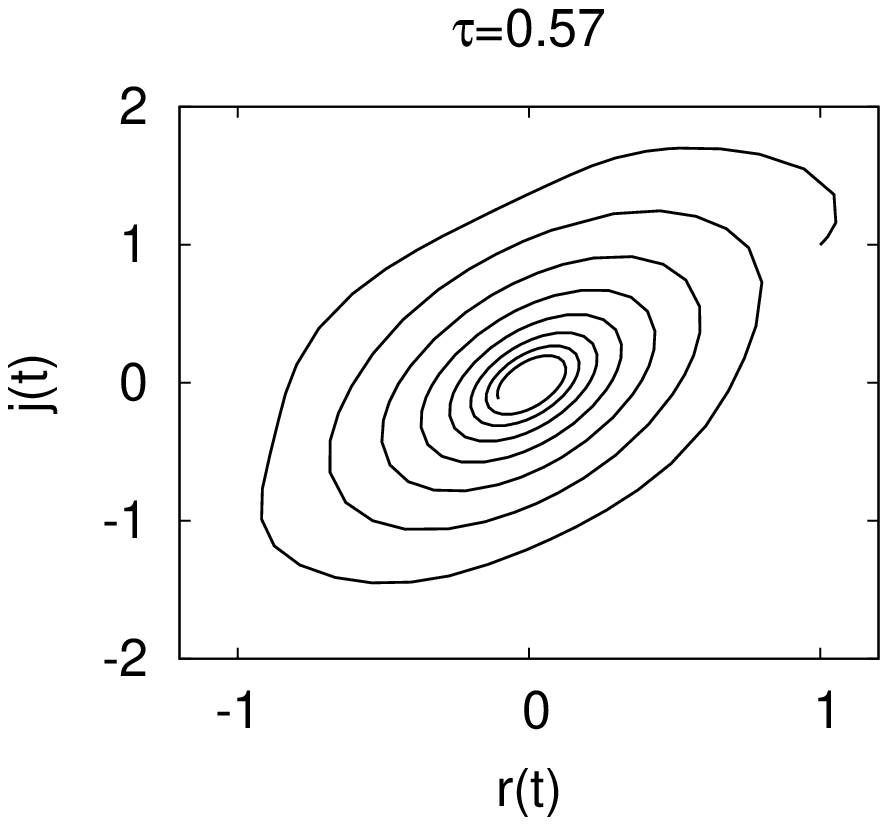}\includegraphics[width=5.5cm]{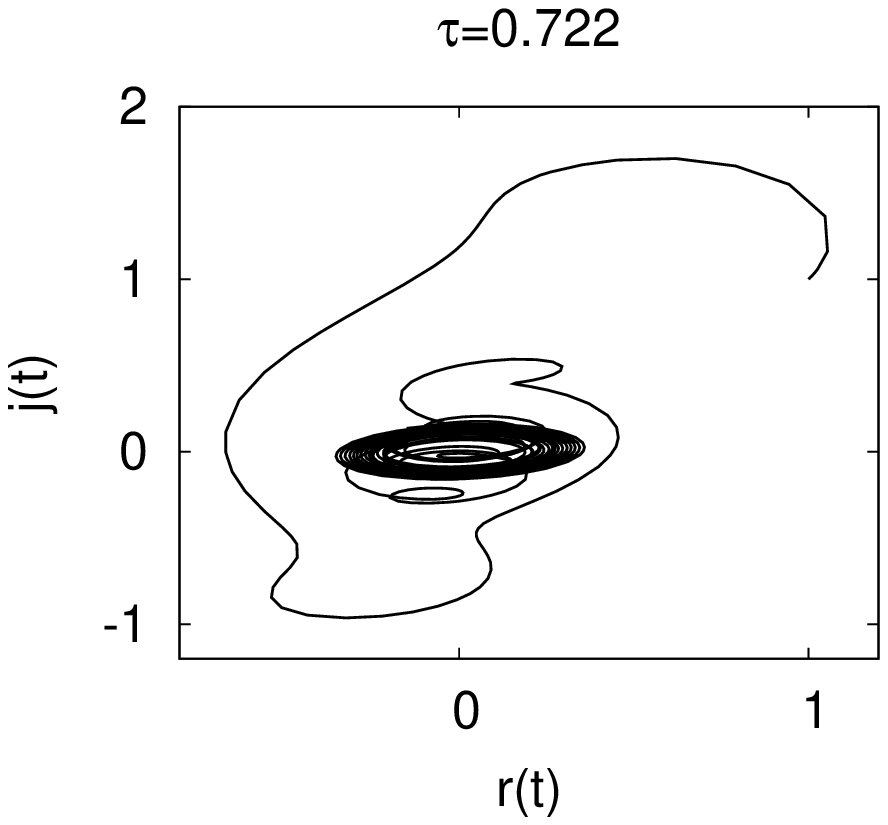}\includegraphics[width=5.5cm]{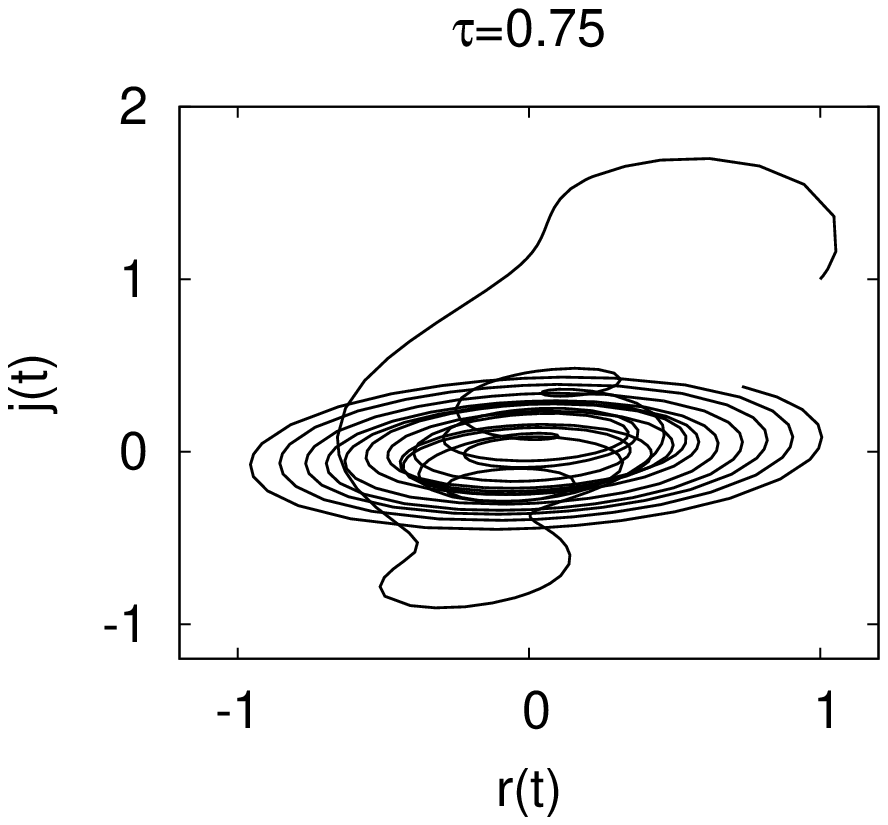}
\caption{Trajectories of Eqs.~(\ref{eq:12}) for $[5,-4,3,-1]$ around two stability switches: the upper three from unstable to stable 
equilibrium, the lower three from stable to unstable one. Middle diagrams indicate the limit cycles}
\label{fig3}
\end{figure}

In Fig.~\ref{fig3} we show the  trajectories around two stability switches, satisfying the assumptions of 
Theorem~\ref{tw2} for the coefficients  $[5,-4,3,-1]$. For delays values below the first switch, the  equilibrium is unstable, 
for delays between the two switches the equilibrium is stable. The left diagrams show the trajectories which converge to the equilibrium. 
The middle one indicate the corresponding limit cycles, and the right ones show the divergent trajectories. 
Similar behavior is obtained e.g. for the reactive parameters $[9,4,-6,1]$, $[7,5,-3,2]$, $[5,-7,4,2]$. Now, at least two of $a_{kl}$ are positive.  

\begin{figure}[ht]
\centering
\includegraphics[width=5.5cm]{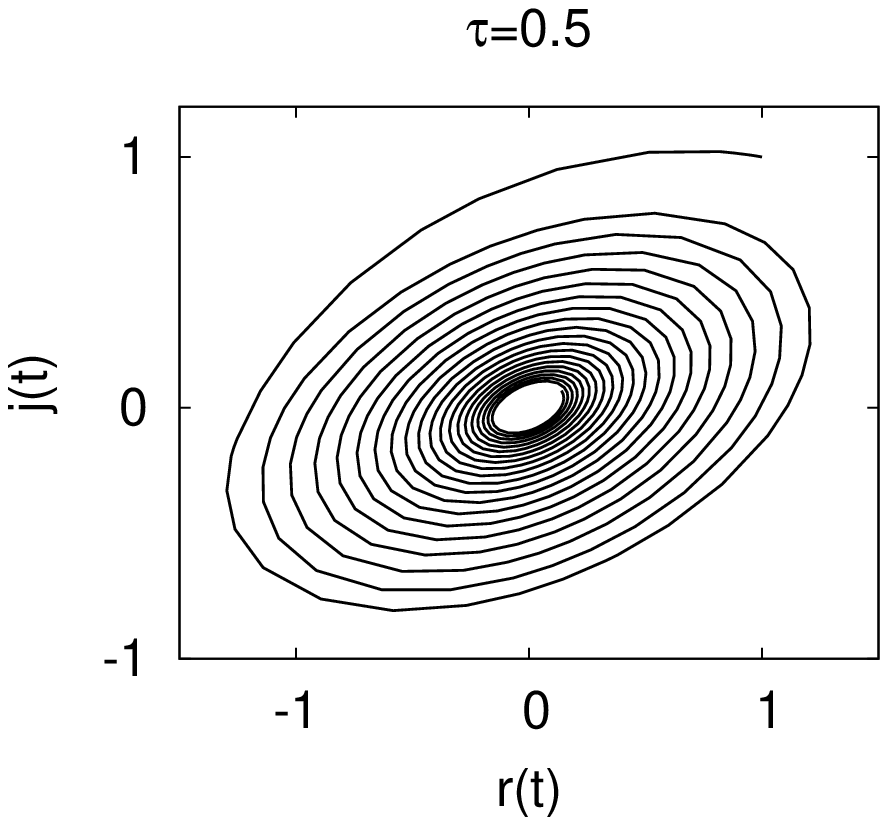}\includegraphics[width=5.5cm]{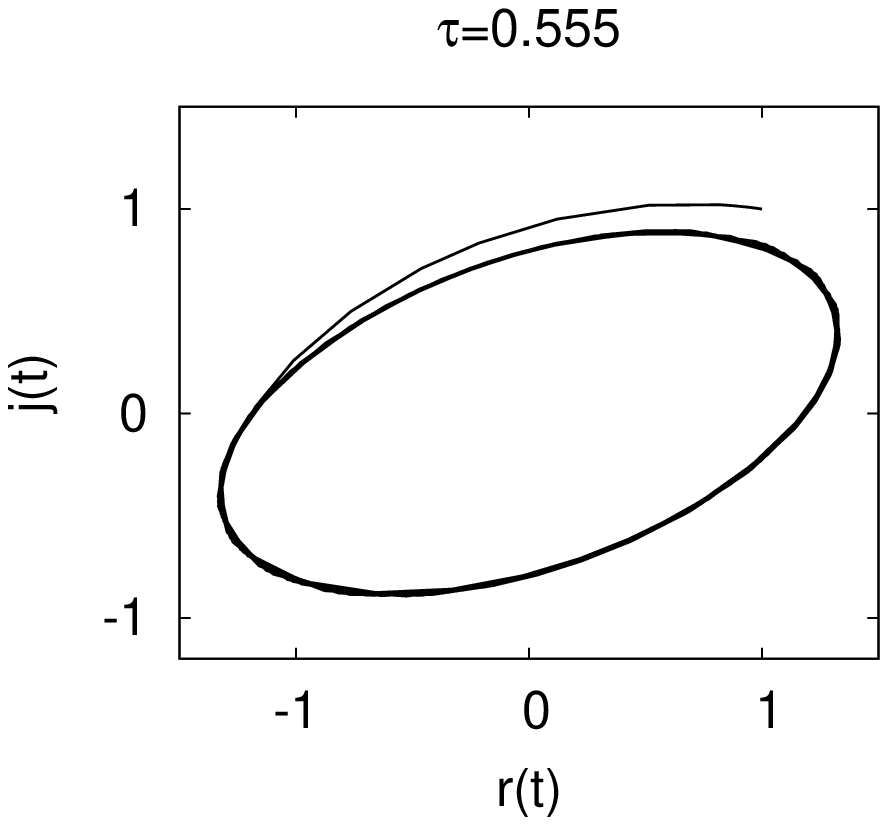}\includegraphics[width=5.5cm]{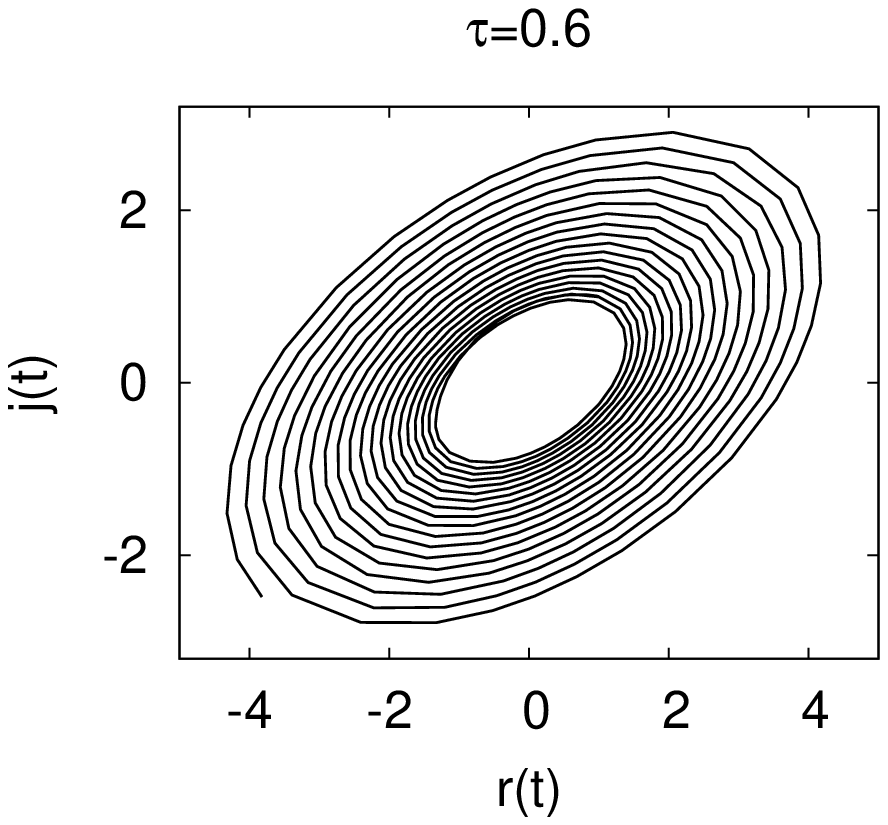}

\centering
\includegraphics[width=5.5cm]{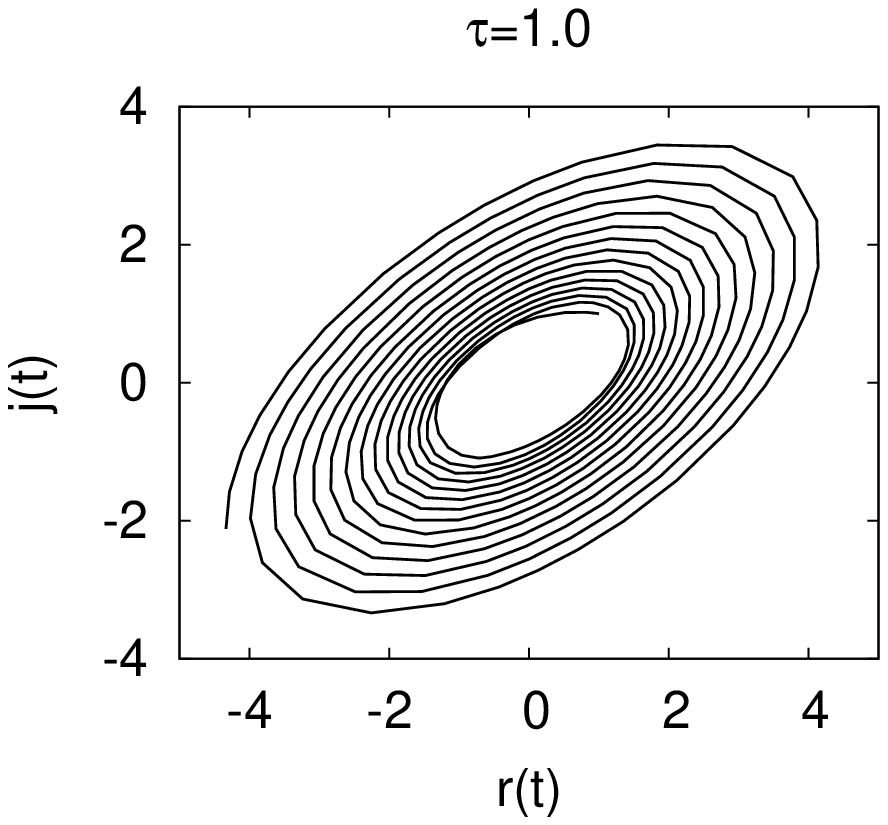}\includegraphics[width=5.5cm]{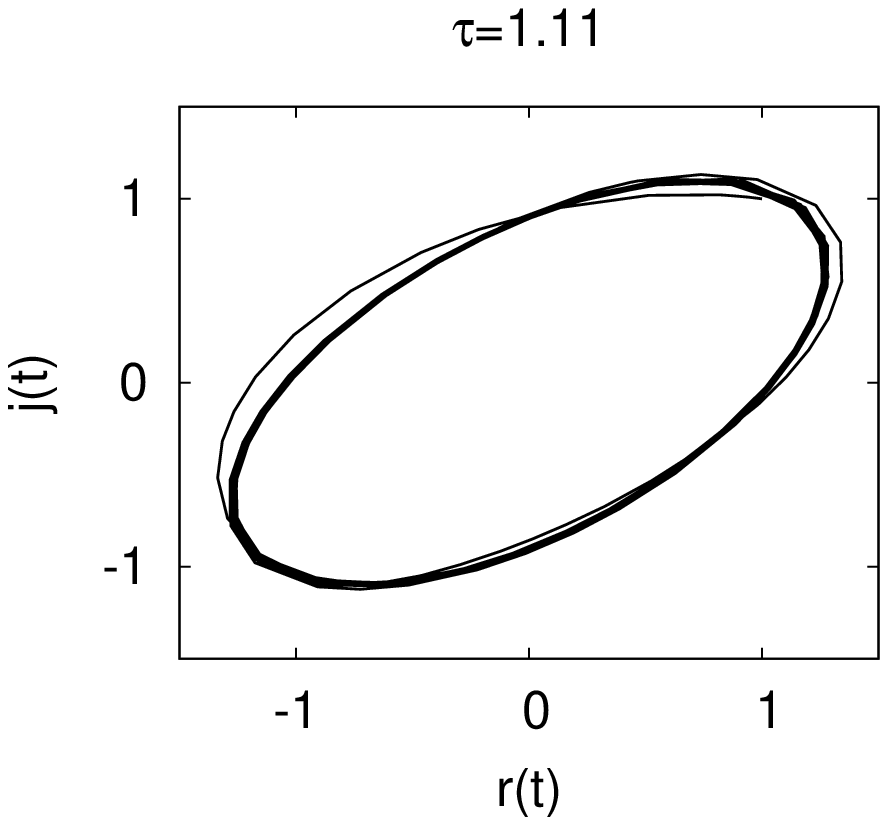}\includegraphics[width=5.5cm]{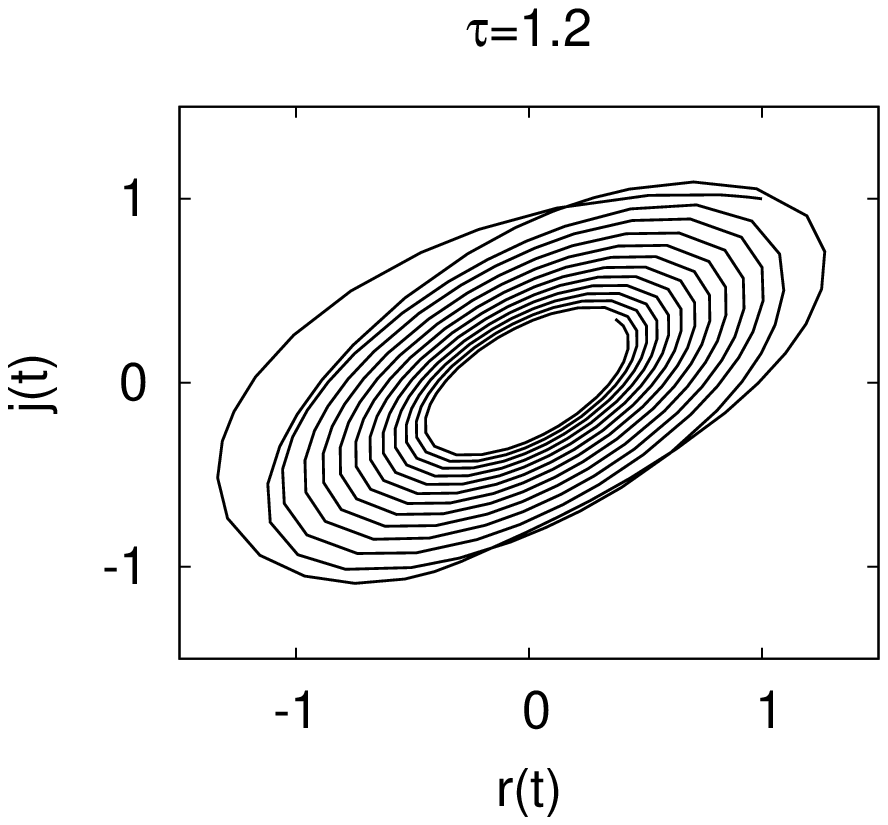}

\centering
\includegraphics[width=5.5cm]{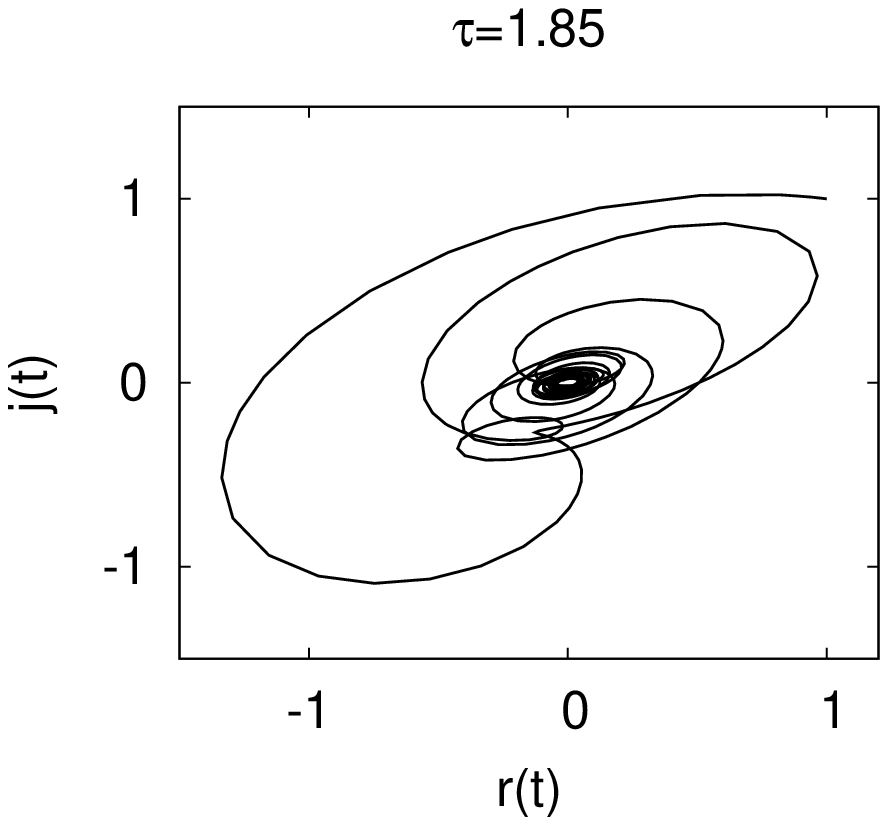}\includegraphics[width=5.5cm]{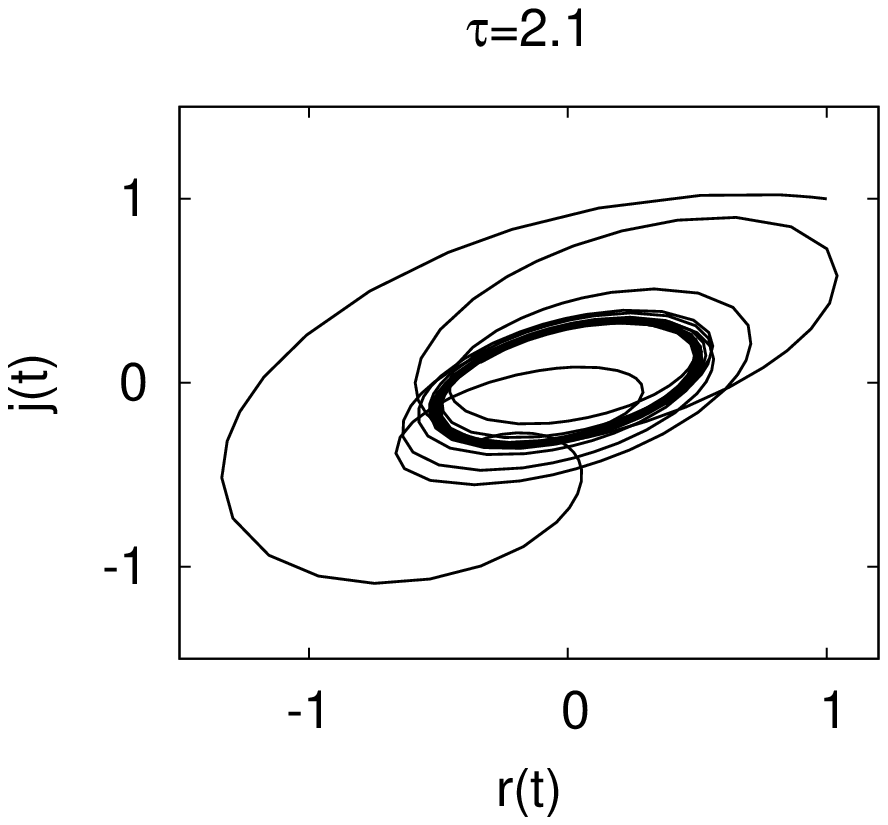}\includegraphics[width=5.5cm]{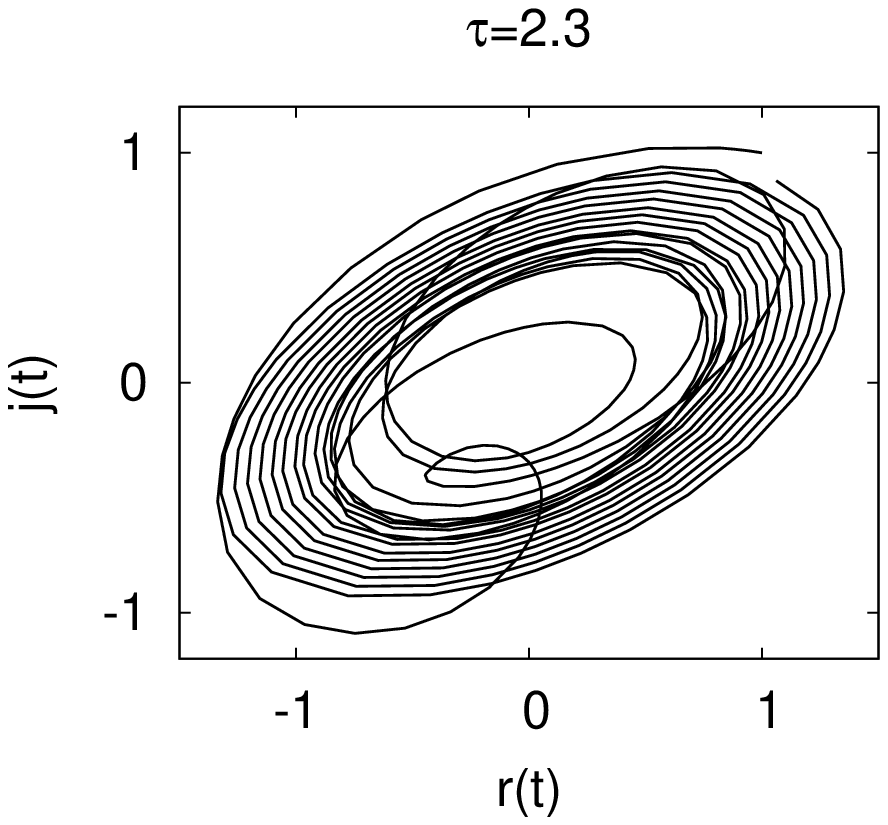}
\caption{ Trajectories of Eqs.~(\ref{eq:12}) for $[-2,-4,3,-2]$ around three stability switches. Middle diagrams indicate the limit cycles }
\label{fig2}
\end{figure}

In Fig.~\ref{fig2} we show the trajectories around three consecutive stability switches 
$\tau \approx 0.555$, $\tau \approx 1.11$, $\tau \approx 2.1$ for the coefficients $[-2,-4,3,-2]$. The left diagram in each row shows the trajectory which converge to equilibrium. The middle one 
approximately visualizes the limit cycle, and the right one shows the diverging trajectories. The corresponding diagrams which present the 
time dependence of $r(t)$, $j(t)$ are similar to those in Fig.~\ref{fig0}, and therefore are not shown. Similar behavior can be obtained e.g. 
for the reactive parameters $[4,6,-5,-5]$, $[5,-7,7,-7]$, $[2,8,-1,-2]$.  Here, at least one of $a_{kl}$ is positive. 

\section{Delay in two terms}
There are four types of models in which two of the four stimuli terms are delayed.   
\begin{equation}
\left\{
\begin{array}{lcccc}
\dot{r}(t) & = & a_{11}r(t - \tau) & + & a_{12}j(t-\tau),\\
\dot{j}(t) & = & a_{21}r(t) & + & a_{22}j(t),
\end{array}
\right. 
\label{eq:27}
\end{equation} 
\begin{equation}
\left\{
\begin{array}{lclcl}
\dot{r}(t) & = & a_{11}r(t - \tau) & + & a_{12}j(t) ,\\
\dot{j}(t) & = & a_{21}r(t - \tau) & + & a_{22}j(t),
\end{array}
\right. 
\label{eq:28}
\end{equation}
\begin{equation}
\left\{
\begin{array}{lcccc}
\dot{r}(t) & = & a_{11}r(t - \tau) & + & a_{12}j(t) ,\\
\dot{j}(t) & = & a_{21}r(t) & + & a_{22}j(t - \tau),
\end{array}
\right. 
\label{eq:31}
\end{equation}
\begin{equation}
\left\{
\begin{array}{lcccc}
\dot{r}(t) & = & a_{11}r(t) & + & a_{12}j(t-\tau) ,\\
\dot{j}(t) & = & a_{21}r(t-\tau) & + & a_{22}j(t). 
\end{array}
\right. 
\label{eq:32}
\end{equation}
The first model corresponds to the situation in which R is reflexive towards both himself and J. In the second model both 
partners consider the R's state with the delay. In the third model both partners are reflexive towards their own states. 
Finally the fourth one corresponds to the situation in which R and J are reflexive towards their partners. 
As before, we can consider symmetric models with the roles of J and R exchanged, but for such models the dynamics is the same. 

We have proved the following   
\begin{tw}\label{tw4}  
\begin{enumerate}
\item Eqs.~(\ref{eq:27}), (\ref{eq:28}) and~\eqref{eq:32} have at most one stability switch.  
\item If $a_{11} + a_{22} = 0$, then Eqs.~(\ref{eq:31}) have at most one stability switch.
\item If $a_{22}<0$, $a_{12}a_{21}<0$ and $a_{11}$ is close enough to $0$, then Eqs.~\eqref{eq:31}  can have arbitrary number of stability switches. 
\end{enumerate}
\end{tw}
For Eqs.~(\ref{eq:31}) with $a_{11} + a_{22} \neq 0$, it is difficult to analyze stability switches with the method presented in Appendix for other systems since the quazi-polynomial has two irreducible trigonometric terms. However, we are able to provide arguments for the statement it does display the same features as Eqs.~\eqref{eq:12} and multiple stability switches are possible. We present these arguments in the proof of Theorem~\ref{tw4} in~Appendix. 

\section{Delay in three and four terms}

There are two families of models in which only one of the stimuli terms is not delayed. 
\begin{equation}
\left\{
\begin{array}{lclcc}
\dot{r}(t) & = & a_{11}r(t - \tau) & + & a_{12}j(t - \tau) ,\\
\dot{j}(t) & = & a_{21}r(t - \tau) & + & a_{22}j(t),
\end{array}
\right. 
\label{eq:40}
\end{equation}
and
\begin{equation}
\left\{
\begin{array}{lcccl}
\dot{r}(t) & = & a_{11}r(t - \tau) & + & a_{12}j(t - \tau) ,\\
\dot{j}(t) & = & a_{21}r(t) & + & a_{22}j(t-\tau). 
\end{array}
\right. 
\label{eq:41}
\end{equation}
In the first case J reacts immediately only on her own state, in the second case only to R's state. 

Finally, there can be the case when both partners react with the delay on both states, that is 
\begin{equation}
\left\{
\begin{array}{lclcl}
\dot{r}(t) & = & a_{11}r(t - \tau) & + & a_{12}j(t - \tau) ,\\
\dot{j}(t) & = & a_{21}r(t-\tau ) & + & a_{22}j(t-\tau). 
\end{array}
\right. 
\label{eq:4_del}
\end{equation}

The quazi-polynomials for systems~\eqref{eq:40} and~\eqref{eq:41} are presented in Appendix. They  contain irreducible trigonometric terms, as in the case of Eqs.~\eqref{eq:31}. This means that we are not able to formulate appropriate theorems, while for Eqs.~\eqref{eq:4_del} we have proved, cf.~Appendix, the following
\begin{tw}\label{tw5}
If system~\eqref{eq:4_del} is unstable for $\tau=0$, then it remains unstable for any $\tau>0$.
If system~\eqref{eq:4_del} is stable for $\tau=0$, then there is an unique switch of stability and the Hopf bifurcation occurs.
\end{tw}
The result of Theorem~\ref{tw5} implies that if all reactions of both partners are equally delayed, then the dynamics of the system is similar to those for one variable model, cf.~e.g.~\cite{Sko}.

Coming back to systems~\eqref{eq:40} and~\eqref{eq:41}, we might think analogously to the previous models and check the dynamics for possible simplifications. For both models all simplifications (cf.~Appendix) give the same results with at most one stability switch. On the other hand, we can provide some arguments that multiple stability switches are possible. We present it in Note to systems with three delays in Appendix.  
However, numerical simulations have not bring any other possibilities than no stability switches or only one stability switch. We do not present such simulations in the paper because they are almost the same as presented for the models with single delay.

\section{Triadic interactions}

One of the possible extensions of the considered models may be obtained by introducing the third person. 
Such generalization has been already suggested in \cite{Str1}, and recently investigated in \cite{Bag}. In particular, 
in the latter work the authors analyzed periodic and quasi periodic solutions for linear and nonlinear models of 
three--person love relations without time delay. 
 
Let Paris (P) enters the relationship of R and J. There are many possible ways of modeling such a triad with delays. We discuss two examples described by the following systems of three DDEs
\begin{equation}
\left\{
\begin{array}{lccclcl}
\dot{r}(t) & = & a_{11}r(t) & + & a_{12}j(t), &&\\
\dot{j}(t) & = & a_{21}r(t-\tau) & + & a_{22}j(t) & + & a_{23}p(t-\tau),\\
\dot{p}(t) & = & a_{32}j(t) & + & a_{33}p(t) , &&
\end{array}
\right.  
\label{eq:301}
\end{equation} 
\begin{equation}
\left\{
\begin{array}{lclcccl}
\dot{r}(t) & = & a_{11}r(t) & + & a_{12}j(t) , &&\\
\dot{j}(t) & = & a_{21}r(t) & + & a_{22}j(t-\tau) & + & a_{23}p(t) ,\\
\dot{p}(t) & = & a_{32}j(t) & + & a_{33}p(t) .&&
\end{array}
\right.  
\label{eq:302}
\end{equation} 
Both models describe the fact that R and P do not interact directly, therefore do not react to the state of the opponent, whereas J needs 
time to reconsider the situation and react. In the model~(\ref{eq:301}) the pragmatic J analyzes and reacts with delay to the 
states of both men, whereas in~(\ref{eq:302}) the reaction of introvertic J to their states is impulsive, immediate, whereas 
the  influence of her states on the speed of change of her actual state is delayed. 

The analysis of the relevant quazi-polynomials which have the third order terms becomes complex.
In a particular case of Eqs.~(\ref{eq:301}) in which R and P react identically to their own states: $a_{11}=a_{33} = a$, we obtain an interesting result. Substituting $s(t):=a_{21}r(t) + a_{23}p(t)$ we achieve
\begin{equation}
\left\{
\begin{array}{lcccc}
\dot{j}(t) & = & a_{22}j(t) & + & 1 \cdotp s(t-\tau) ,\\
\dot{s}(t) & = & (a_{12}a_{21}+ a_{23}a_{32})j(t) & + & as(t) .
\end{array}
\right.
\label{eq:303}
\end{equation} 
Thus, we can perceive this situation as if J would interact with an ''averaged'' man whose state is $s$, and is influenced by this state with the delay. We obtain the model~(\ref{eq:23}) with the roles of R and J interchanged, in which only one stability switch is possible. Thus, when  both men have the same attitude 
towards themselves, there is at most one stability switch. 

Below we present an example of numerical results for the second model~\eqref{eq:303} for the case with three stability switches. 
\begin{figure}[ht]
\centering
\includegraphics[width=5.5cm]{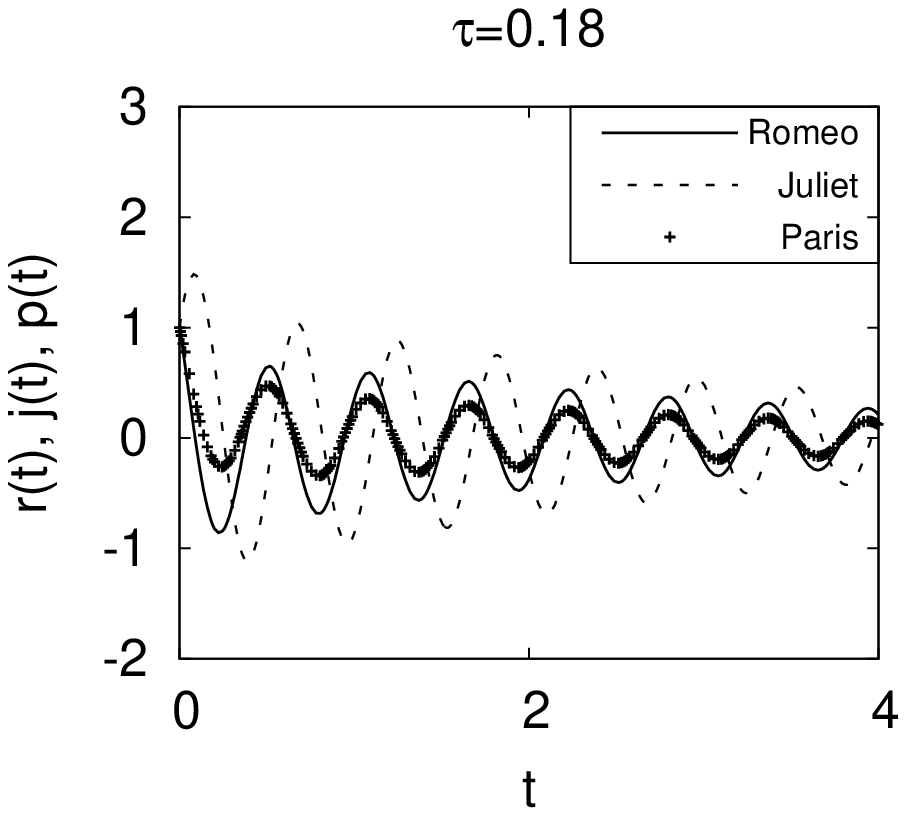}\includegraphics[width=5.5cm]{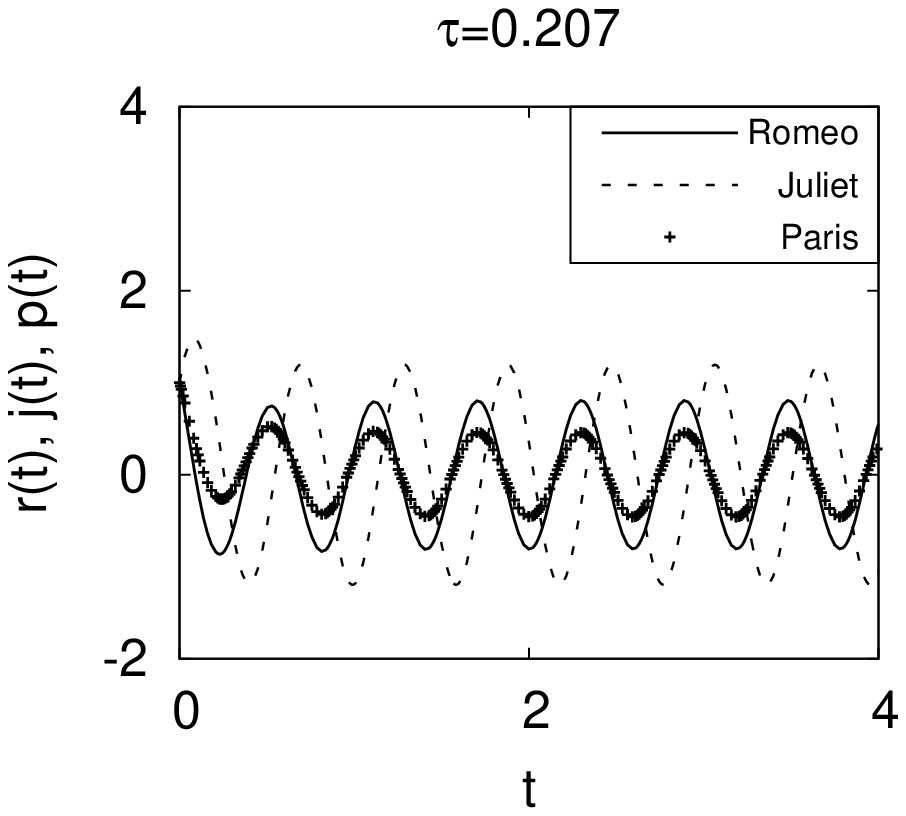}\includegraphics[width=5.5cm]{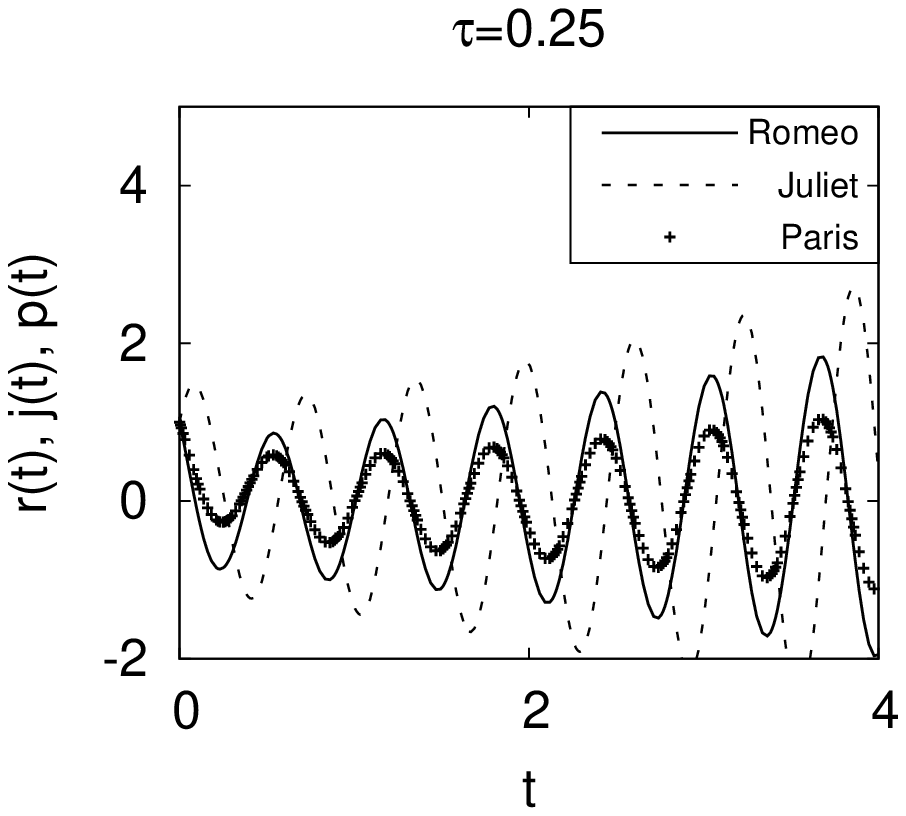}

\centering
\includegraphics[width=5.5cm]{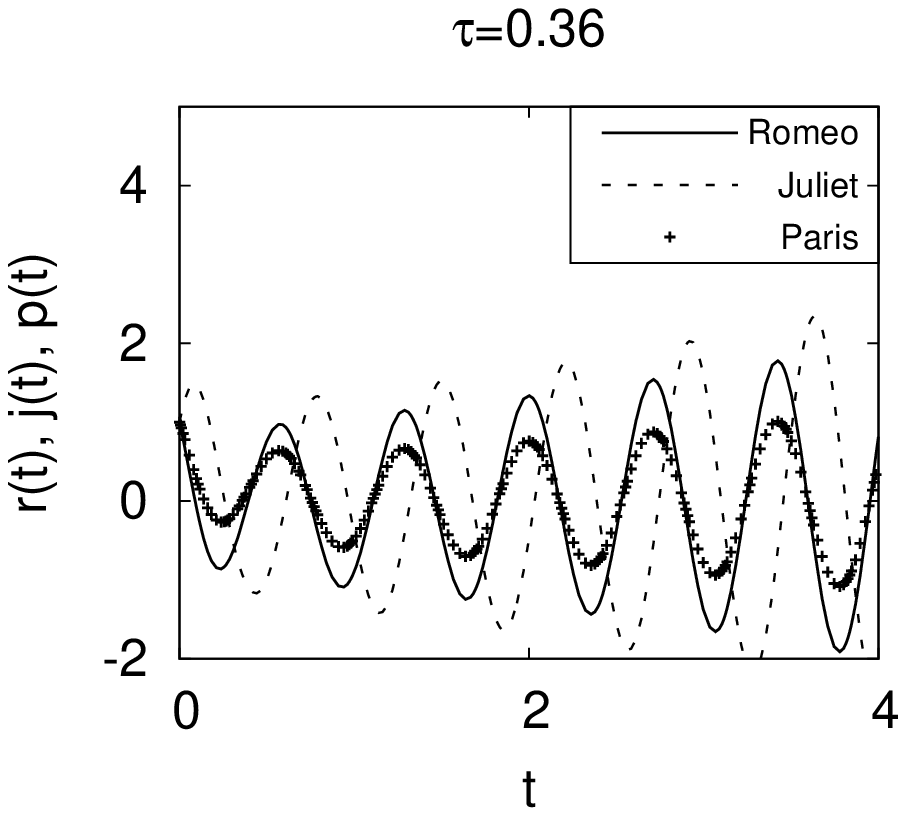}\includegraphics[width=5.5cm]{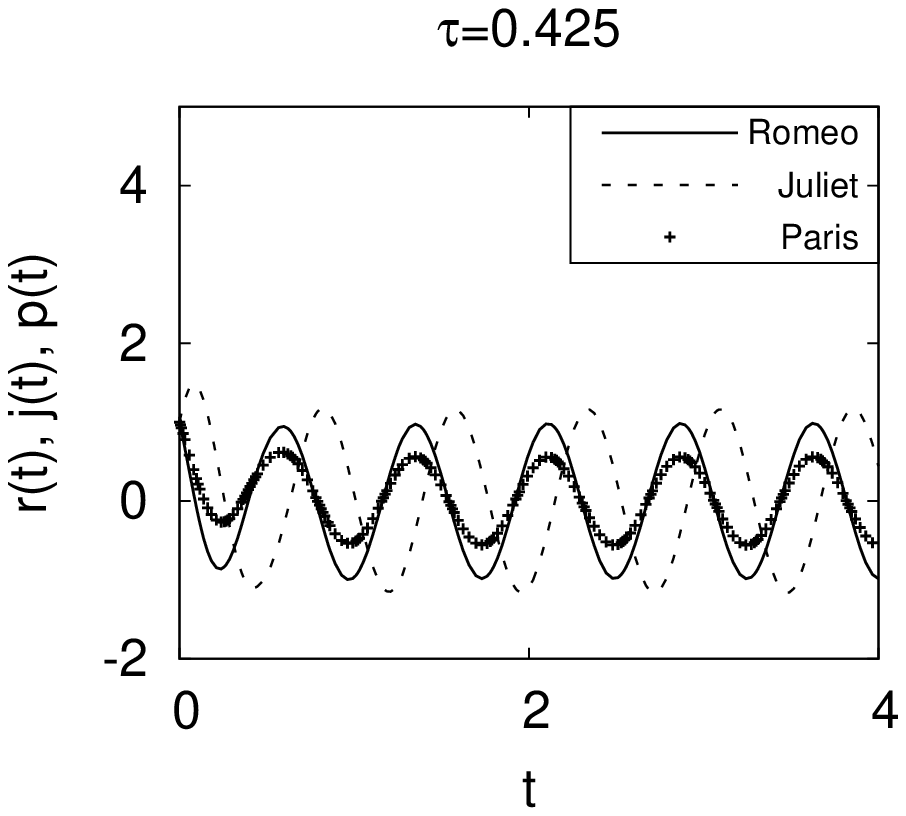}\includegraphics[width=5.5cm]{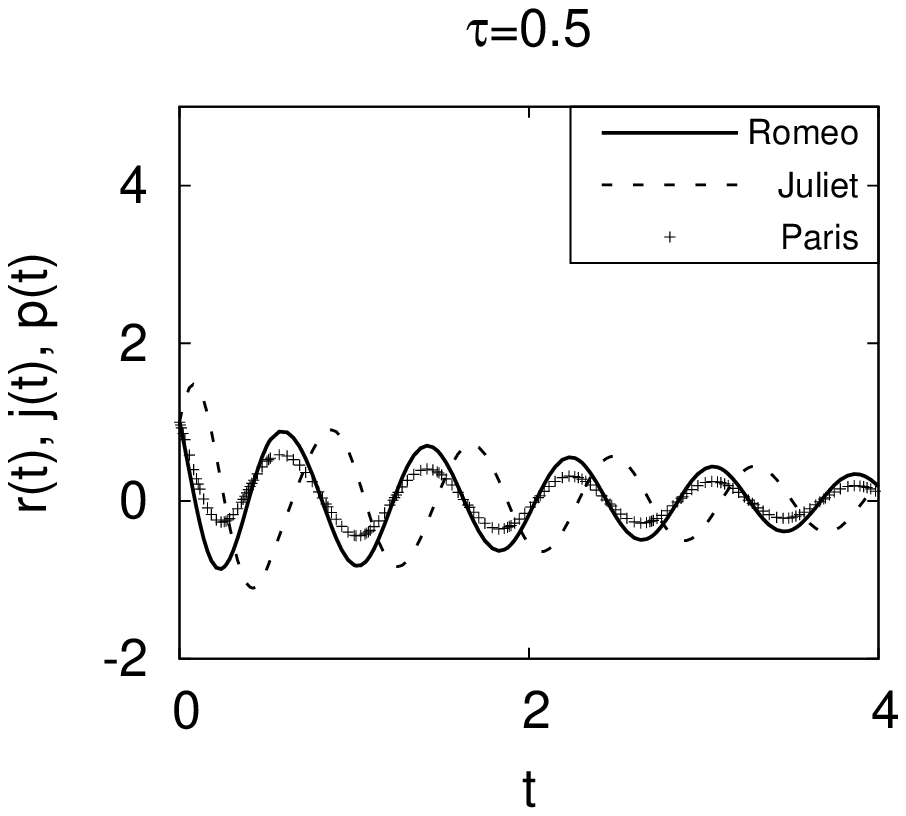}

\centering
\includegraphics[width=5.5cm]{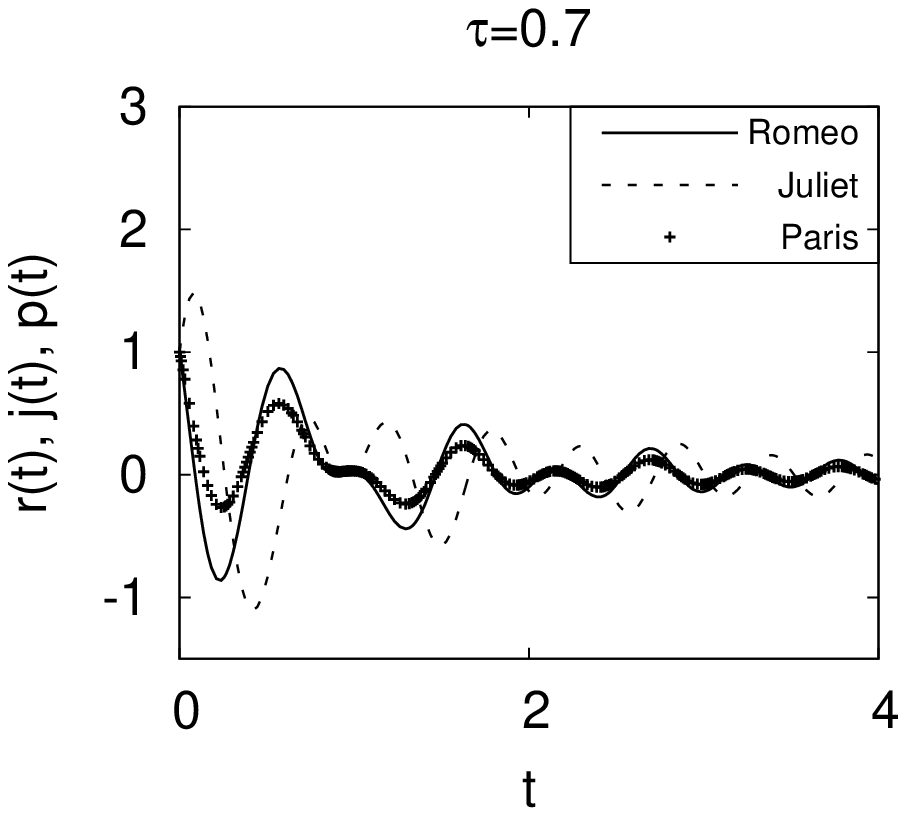}\includegraphics[width=5.5cm]{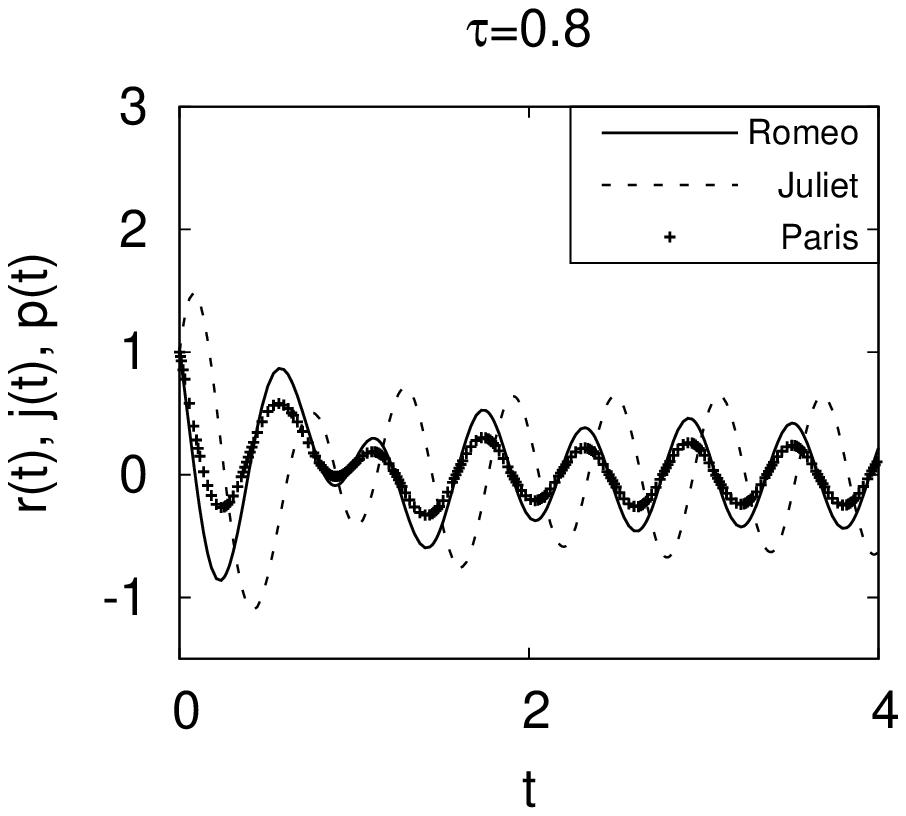}\includegraphics[width=5.5cm]{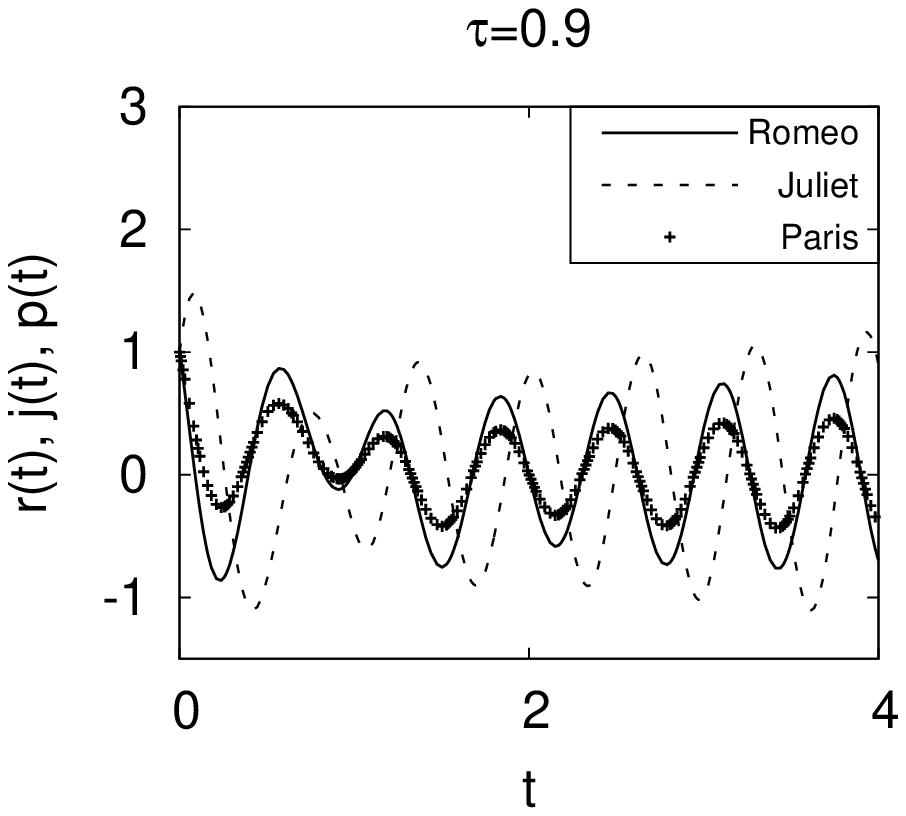}
\caption{Solutions to Eqs.~(\ref{eq:302}) for 
$[a_{11},a_{12},a_{21},a_{22},a_{23},a_{32},a_{33}]=[-28,-74,76,-35,76,-42,-28]$}
\label{figtriada1}
\end{figure}
Analyzing the relevant quazi-polynomials we can take advantage from the fact that stability of the system depends continuously on parameters in order to find 
examples of triads with multiple stability switches. For example, in order to obtain a system with $n$ stability switches, 
we may consider a subsystem with only $a_{11}, a_{12}, a_{21}, a_{22}$ to be non-zero parameters and obtain R--J system with $N$ 
stability switches. Then, we may put $a_{32}=a_{12}$, $a_{33}=a_{11}$, $a_{23}=a_{21}$ to create a system to be symmetrical towards 
R--P axis. This system will share the number of switches with the subsystem. From this point we may try to change values of chosen parameters until the number of stability switches jumps to another value. It can be shown that multiple stability switches can also occur for adequately chosen intervals of the parameters of the models. 
\vskip 0.1cm
We have used this method to obtain a system with three stability switches for parameters $a_{11}=-28$, $a_{12}=-74$, $a_{21}=76$, $a_{22}=-35$, $a_{23}=76$, $a_{32}=-42$, $a_{33}=-28$. For these parameters, the model~(\ref{eq:302}) describes a situation when after deliberation J reacts negatively to her own states, but is quick in response to both R and P, and follows their emotions. R and P are similar in reactions and tend to neutralize the relationship: the more satisfied J is, the more satisfaction they lose and the opposite. They are the same ambivalent towards themselves.  
The behavior of this system is shown in Fig.~\ref{figtriada1}.

\section{Summary and conclusions}
In the paper we have studied a general family of linear models of dyadic relationships with one constant time delay 
present on the rhs of the systems. We have used the mathematical approach presented in~\cite{Coo}. 
The only linear systems of two DDEs which can not be completely analyzed using the presented method are those for which the quazi-polynomial has the irreducible term  
 $a \text{e}^{-\lambda\tau} + b\text{e}^{-2\lambda\tau}$, $a^2 + b^2 \ne 0$, which appears if the delays is present in the 'cross' terms of the rhs
of the considered system (since then the auxiliary function $F(y)$ contains trigonometric terms, and in general analytic considerations are very difficult). 

The results presented in the paper, stated mainly in Theorems~\ref{tw1}-\ref{tw4}, show that the joint strength of reactions of partners (where the joint strength means here the product of appropriate strengths of both partners) 
to the partner's state has greater impact on the system dynamics than the joint strength of reactions on the partner's own state. This is reflected in the fact that for $|a_{11}a_{22}|<|a_{12}a_{21}|$, that is for the joint strength of reactions of partner's own state weaker than the joint strength of reactions to partner's state, the dynamics is typically much simpler. More precisely, Theorem~\ref{tw1} for $|a_{11}a_{22}|>  |a_{12}a_{21}|$ implies that the system~\eqref{eq:12} can have at most one stability switch, while Theorem~\ref{tw2} for the same system with opposite inequality can exhibit multiple stability switches. Similarly, for Eq.~\eqref{eq:23} Theorem~\ref{tw3} yields no stability switches for $|a_{11}a_{22}|>|a_{12}a_{21}|$ and possibility of one stability switch for the opposite inequality. 

The fact that stability of the system depends continuously on the parameters leads to a general observation that 
adding a term 
with the interaction parameter close enough to 0 to a given system does not change its stability properties. 
However, there is a substantial difference between systems depending on which term is delayed. 
Among all the analyzed two--person interactions, only those described by the system~\eqref{eq:12} provide 
the possibility of multiple stability switches. In other words, multiple stability switches are possible 
exclusively when the stimulus of a player to his own state is delayed, and the remaining three stimuli are not. 
Otherwise, i.e. when only the stimulus to the state of the partner  is delayed, or when at least two stimuli are 
delayed, then at most one stability switch occurs. 

In this paper we have considered the interactions modeled by linear systems. It can be of interest to consider nonlinear 
interactions between the partners.  
On the other hand, one could consider the most general linear model with arbitrary discrete delays. Let 
$a_k$, $b_k$, $c_k$, $d_k$ denote the constant reactive coefficients, with for example 
$k=1,...,n_{rr}$ for $a_k$ etc., $\tau_k^{rr}\ge 0$, $\tau_k^{rj}\ge 0$,  $\tau_k^{jr}\ge 0$,   $\tau_k^{jj}\ge 0$ denote the 
relevant delays, corresponding to different interaction terms in the relationship R--J. The general model reads
$$
\frac{d}{dt} r(t) = \sum_{k=1}^{n_{rr}} a_k r(t-\tau_k^{rr}) + \sum_{k=1}^{n_{rj}} b_k j(t-\tau_k^{rj}) , \\
$$
$$
\frac{d}{dt} j(t) = \sum_{k=1}^{n_{jr}} c_k r(t-\tau_k^{jr}) + \sum_{k=1}^{n_{jj}} d_k j(t-\tau_k^{jj}) ,
$$
and covers all the situations considered in this paper. In particular one could consider the situation when 
the (absolute) values of the reactive coefficients decrease with increasing values of the delays, which would mean that the 
states most remote in the past influence less the actual dynamics of the behavior of the partners.

\section{Appendix}
In the proofs presented below we use the approach developed in~\cite{Coo}, cf.~also~\cite{Sko}. If $W$ is the quazi-polynomial for the system of  DDEs with a single delay $\tau$, then $W(\lambda)=P(\lambda)+Q(\lambda)\text{e}^{-\lambda \tau}$, where $P$ and $Q$ are polynomials. 
If all zeros of $W$ (that is characteristic values or eigenvalues of the system) are on the left-hand side of the complex plane, then the system is stable. If at least one zero is on the right-hand side of the complex plane, then the system is unstable. Because zeros of $W$ depend continuously on the model parameters, the change of stability is possible only when some zero crosses imaginary axis either from the left to the right or from the right to the left-hand side of the complex plane. This is the reason we are looking for purely imaginary zeros of $W$.
Substituting $\lambda:=i\omega$ one gets
\[
W(i\omega)=P(i\omega)+Q(i\omega)\text{e}^{-i\omega\tau} .
\]
Looking for purely imaginary eigenvalues we solve $W(i\omega)=0$, that is
$P(i\omega)=-Q(i\omega)\text{e}^{-i\omega \tau}$. Taking the norm of both sides 
and substituting $\omega^2=y$ we obtain the auxiliary function $F$  defined as
\begin{equation}\label{F}
F(y) :=  |P(i\sqrt{y})|^2-|Q(i\sqrt{y})|^2 ,  
\end{equation} 
and we see that positive zeros of $F$ give the eigenvalues of the studied  system of DDEs. 
We also know that the derivative of $F$ at its zero point determines the direction of movement of the corresponding eigenvalues in the complex plane, cf.~\cite{Coo, Sko}. In fact, if this derivative is positive, then the eigenvalue crosses the imaginary axis from the left to the right-hand side of the complex plane, while for the negative derivative the movement is in the opposite direction.

\medskip 

\noindent
{\bf{Proof of Theorem~\ref{tw1}.}} 
The quazi-polynomial for Eqs.~\eqref{eq:12}  has the form
\[
W(\lambda) = \lambda^2 - a_{22}\lambda -	 a_{12}a_{21} + \text{e}^{-\lambda\tau}(a_{11}a_{22} - \lambda a_{11}) .
\]
The auxiliary function $F$ reads
\begin{equation}\label{F12}
F(y) =  y^2 + y(a_{22}^2 - a_{11}^2 + 2a_{12}a_{21}) + a_{12}^2a_{21}^2 - a_{11}^2a_{22}^2 .  
\end{equation}
If the system is unstable for $\tau=0$, then the corresponding characteristic polynomial, which coincides with $W(\lambda)$ for $\tau=0$, has zeros 
on the right-hand side of the complex plane. As $\tau$ increases, roots of the quazi-polynomial 
should belong to the left-hand side of the complex plane in order to stabilise the system. 

Under the assumption $|a_{11}a_{22}|>|a_{12}a_{21}|$ the function $F$ has only one positive zero $y_0$ and its derivative 
$$F'(y_0) = 2y_0 + a_{22}^2 – a_{11}^2 + 2a_{12}a_{21} > 0 ,$$
which means that the transition from the right to 
the left-hand side is not possible, cf.~\cite{Sko}.   

This implies that if the system is unstable for $\tau=0$, 
then it remains unstable for all $\tau>0$. Moreover, if the eigenvalue enters the right-hand complex half-plane, it 
 remains there forever, yielding that only one stability switch is possible. 

If the system is stable for $\tau=0$, 
i.e.  $a_{11}a_{22} - a_{12}a_{21} > 0$ and $a_{11} + a_{22} < 0$ under the assumption of Theorem~\ref{tw1}, 
then all characteristic values are on the left-hand complex half-plane and can cross to the right-hand complex half-plane which implies destabilisation. There exists a critical 
$\tau_c$ for which HB occurs and the system is unstable for larger delays.  
More precisely, let $\omega_0 = \sqrt{y_0}$, where $y_0$ is the unique positive zero of $F$. Then the HB occurs for the 
smallest $\tau$  for which $W(i\omega_0) = 0$, i.e. for the smallest $\tau$ such that  
\begin{equation}\label{re-im}
\left\{
\begin{array}{lcccl}
\Re(W(i\omega_0))  & = & -\omega_0^2 -a_{12}a_{21} - a_{11}\omega_0 \sin(\omega_0 \tau) + a_{11}a_{22}\cos (\omega_0 \tau) & = & 0 ,\\
\Im(W(i\omega_0)) & =  & -a_{22}\omega_0 - a_{11}a_{22}\sin(\omega_0 \tau) - a_{11}\omega_0 \cos(\omega_0 \tau) & = & 0,
\end{array}
\right. 
\end{equation}
i.e. for $\tau_c$ determined by the values of sine and cosine 
\[
\cos(\tau_c\omega_0)=\frac{a_{12}a_{21}a_{22}}{a_{11}(a_{22}^2+\omega_0^2)}, \quad 
\sin (\tau_c\omega_0)=- \frac{\omega_0(\omega_0^2 +a_{22}^2-a_{12}a_{21})}{a_{11}(a_{22}^2+\omega_0^2)} ,
\]
which implies   
\[
\tau_c = 
\left\{ \begin{array}{ccc}
\frac{1}{\omega_0} \arccos a_0 & \text{for} & a_{11}<0 ,\\
2\pi - \frac{1}{\omega_0} \arccos a_0 & \text{for} & a_{11}>0  .
\end{array}
\right.
\]
$\Box$

\medskip

\noindent
{\bf{Proof of Theorem~\ref{tw2}.}} 
If $|a_{12}a_{21}|>|a_{11}a_{22}|$, then stability switches can occur only if the auxiliary function $F$ defined by Eq.~\eqref{F12}  has two positive roots 
$0<y_0 < y_1$, i.e. if
\[
a_{22}^2 - a_{11}^2 + 2a_{12}a_{21}>0 ,\ \
(a_{22}^2 + 2a_{12}a_{21} - a_{11}^2)^2 > 
4\left((a_{12}a_{21})^2 - (a_{11}a_{22})^2\right) .
\]
Since $F'(y) = 2y + (a_{22}^2 - a_{11}^2 + 2a_{12}a_{21})$ we have 
\[
F'(y_0) = -\sqrt{\Delta} < 0\quad \text{and} \quad
F'(y_1) = +\sqrt{\Delta} > 0,
\]
where $\Delta$ is the discriminant of $F$,
and therefore we deduce that for $\omega_0$ the eigenvalue can move from the right half-plane to the left one, while for $\omega_1$ -- in the opposite direction. 

Eqs.~\eqref{re-im} yield 
\begin{equation}
\left\{
\begin{array}{lcc}
\sin(\omega_j \tau_j) & = & \frac{-\omega_j(\omega_j^2 + a_{22}^2 + a_{12}a_{21})}{a_{11}(\omega_j^2 + a_{22}^2)} , \\
\cos(\omega_j\tau_j) & = & a_j ,   
\end{array}
\right.
\label{eq:19}
\end{equation}
with $a_j = \frac{a_{12}a_{21}a_{22}}{a_{11}(\omega_j^2 + a_{22}^2)}$
for $ j=0, 1$. For arbitrary $j$ let $\tau_{j0}>0$ be the first critical delay for which Eqs.~\eqref{eq:19} are   satisfied. 

We obtain two sequences of critical delays for which eigenvalues can cross the imaginary axis 
\begin{equation}
\tau_{jn} = \tau_{j0} + \frac{2n\pi}{\omega_j}, \quad n \in \mathbb{N}, \quad j=0, \; 1.
\label{eq:20}
\end{equation}
There are two possibilities
\begin{enumerate}
\item If the system in unstable for $\tau=0$, then at least one zero is on the right-hand side of the imaginary 
axis. Stability switch occurs if the eigenvalue crosses the imaginary axis from the right to the left-hand side complex plane through the point 
$i\omega_0$. This can happen only for the sequence $(\tau_{j0})$.  Thus,
\[
\tau_{00} < \tau_{10}. 
\]
In this case there are at least two switches, because for sufficiently large $\tau$ the system is unstable, cf.~\cite{Coo, Sko}. 
In order to have more switches, the system which is unstable for $\tau=0$ should satisfy 
\[
\tau_{00} < \tau_{10} < \tau_{01},
\]
which gives two corresponding inequalities in (\ref{eq:6000}).
\item If the system is stable for $\tau=0$, then it has all eigenvalues in the lhs of the complex plane. 
If $\tau_{10} < \tau_{00}$, then at least one stability switch occurs. In order to have three switches the 
following (together with the assumptions listed before) conditions should be satisfied
\[
\tau_{10}<\tau_{00}<\tau_{11} .
\]
\end{enumerate}
Analogously one can obtain the set of sufficient conditions for 
arbitrary number $N$ of stability switches.

\noindent
$\Box$

Now, consider the possibility of $5$ stability switches. Let the 
system be stable for $\tau=0$, then there are $5$ switches if 
\[
\left\{
\begin{array}{l}
\Delta > 0,\\
|a_{12}a_{21}|>|a_{11}a_{22}|,\\
a_{11} + a_{22} <0 \ \wedge \  a_{11}a_{22} - a_{12}a_{21}>0,\\
a_{11}^2 - a_{22}^2 - 2a_{12}a_{21} > 0,\\
{\omega_0}\arccos a_1 < {\omega_1}\arccos a_0 < {\omega_0}\arccos a_1 + {2\pi}{\omega_0},\\
{\omega_0}\arccos a_1 +2\pi\omega_0 < {\omega_1}\arccos a_0 + 2\pi\omega_1 < {\omega_0}\arccos a_1 + {4\pi}{\omega_0},
\end{array}
\right.
\]
Five stability switches occur for example for the system~(\ref{eq:12}) with the coefficients $[-1, 3, -2, 1]$. 
The corresponding approximate values of the stability switches are 
$$\tau_1\approx 0.92, \ \tau_2\approx1.4, \  \tau_3\approx 3.3, \  \tau_4\approx 4.2, \ \tau_5\approx 5.46. $$ 
The system is stable for $\tau \in [0, \tau_1) \cup (\tau_2, \tau_3) \cup (\tau_4, \tau_5).$

\medskip  

\noindent
{\bf{Note to Remark~1.}}
The quazi-polynomial for Eqs.~(\ref{eq:410}) reads
\begin{equation}
W(\lambda) = \lambda^2 + (-a_{22} - a_{13})\lambda + (-a_{11}\lambda + a_{11}a_{22})\text{e}^{-\lambda \tau} + a_{13}a_{22} - a_{12}a_{21},
\label{eq:43}
\end{equation}
and the corresponding auxiliary function $F(y) = y^2 + (a_{22}^2 + a_{13}^2 +2a_{12}a_{21} - a_{11}^2)y + (a_{13}a_{22} - a_{12}a_{21})^2 - (a_{11}a_{22})^2$. 
We see that for $a_{13}=0$ both functions are the same as for Eqs.~\eqref{eq:12}. 
This yields that for appropriate values of the model parameters there can be arbitrary number of stability switches, as for Eqs.~\eqref{eq:12}.

\medskip

\noindent
{\bf{Proof of Theorem~\ref{tw3}.}}
The quazi-polynomial and the auxiliary function $F$ have the form
\[
W(\lambda) = \lambda^2 - (a_{11} + a_{22})\lambda + a_{11}a_{22} - a_{12}a_{21}\text{e}^{-\lambda\tau}, \]
\[  F(y) = y^2 + (a_{11}^2 + a_{22}^2)y + (a_{11}a_{22})^2 - (a_{12}a_{21})^2. 
\]
If $|a_{12}a_{21}| > |a_{11}a_{22}|$, then $F(0)<0$ and there is one positive zero $y_0$, for which $F'(y_0) = \sqrt{\Delta} > 0$, therefore eigenvalues can cross the imaginary axis only from the left to the right. Thus, 
\begin{enumerate}
\item if the system is unstable for  $\tau = 0$, then te same holds for all $\tau\geq 0$;
\item if the system is stable for $\tau = 0$, then there exists $\tau_c$, for which HB occurs, and for $\tau > \tau_c$ the system 
becomes unstable (the critical value $\tau_{cr}$ can be calculated as in the proof of Theorem~\ref{tw1}).
\end{enumerate}
If $ |a_{12}a_{21}| < |a_{11}a_{22}|$, then the stability switches are possible only if $F(y)$ has two positive roots. 
The necessary condition $a_{11}^2 + a_{22}^2 < 0$ for this situation excludes such possibility. 

\noindent
$\Box$

\medskip 

\noindent
{\bf{Proof of Theorem~\ref{tw4}.}}
\begin{enumerate}
\item For the models~(\ref{eq:27}) and~(\ref{eq:28}), the quazi-polynomial is of the same form 
\[
W(\lambda) = \lambda^2 - a_{22}\lambda + (-a_{11}\lambda + a_{11}a_{22} - a_{12}a_{21})\lambda \text{e}^{-\lambda 
\tau} ,\]and the corresponding auxiliary function 
$F(y) = y^2 + (a_{22}^2 - a_{11}^2)y - (a_{11}a_{22} - a_{12}a_{21})^2$. If $a_{11}a_{22} - a_{12}a_{21} \neq 0$, 
then the free term of the function $F$ is always negative and it has one positive root, which implies there is at most 
one stability switch. 
\item For the model~(\ref{eq:31}) the quazi-polynomial reads  
\begin{equation}\label{pol:31}
W(\lambda) = \lambda^2 - a_{12}a_{21} + (-a_{11} - a_{22})\lambda \text{e}^{-\lambda \tau} + a_{11}a_{22}\text{e}^{-2\lambda\tau}  .
\end{equation}
If in this particular relationship, partners are focused on their own states from the past but with the opposite result, which means $a_{11} + a_{22} = 0$, the quazi-polynomial~\eqref{pol:31} reduces to
\[
W(\lambda) = \lambda^2 - a_{12}a_{21} + a_{11}a_{22}\text{e}^{-2\lambda\tau} 
\]
which enables the analysis as above. The relevant auxiliary function is of the form 
$F(y) = y^2 - 2a_{12}a_{21}y + (a_{12}a_{21})^2 - (a_{11}a_{22})^2$, and it has two roots for all the parameters, 
since $\Delta(F) = 4(a_{11}a_{22})^2$ is always positive. Hence,

a. if $a_{12}a_{21} < a_{11}a_{22}$, the free term of the function $F$ is negative, which implies the function has 
a single positive root, so that there is no more than a single switch possible; 

b. if $a_{12}a_{21} > a_{11}a_{22}$, the system never starts from stability, because the conditions imply jointly 
that $a_{11}a_{22} - a_{12}a_{21} < 0$. Thus, if we look for more than one stability switch, 
the condition $\tau_{00} < \tau_{10} < \tau_{01}$ must be satisfied, which implies $\frac{\pi}{2\omega_0} < \frac{\pi}{2\omega_1} < \frac{\pi}{\omega_0}$, where $\omega_0 = (a_{12}a_{21} - a_{11}a_{22})^{\frac{1}{2}}$, $\omega_1 = (a_{12}a_{21} + a_{11}a_{22})^\frac{1}{2}$, but $0 < \omega_0 < \omega_1$ so that the first inequality is never satisfied. 

On the other hand, if $a_{11}=0$, then quazi-polynomial~\eqref{pol:31} reduces to
\[
W(\lambda) = \lambda^2 - a_{12}a_{21}  - a_{22}\lambda \text{e}^{-\lambda \tau} ,
\] 
which enables the analysis, as before. The auxiliary function 
$F(y)=y^2-(2a_{12}a_{21}+a_{22}^2)y+a_{12}^2a_{21}^2$ has two positive roots under the assumption 
$a_{22}^2>4|a_{12}||a_{21}|$ and then multiple stability switches can occur, as in Theorem~\ref{tw4}. 
Due to the continuous dependence on the parameters multiple switches occur also for $a_{11}$ sufficiently close to $0$.  
\item  For the model~(\ref{eq:32}) there is
\[
W(\lambda) = 
\lambda^2  - (a_{11} + a_{22})\lambda  + a_{11}a_{22} - \text{e}^{-2 \lambda \tau}a_{12}a_{21} .
\]
Here we have $F(y) = y^2 + (a_{11}^2 + a_{22}^2)y + (a_{11}a_{22})^2 – (a_{12}a_{21})^2$ and respectively $\Delta(F) = (a_{11}^2–a_{22}^2)^2 + 4(a_{12}a_{21})^2 > 0 $ for any parameters which means $F$ has two roots. Note that if we expect both the roots to be positive, their sum shall be positive, while sum of the roots equals $-(a_{11}^2 + a_{22}^2)$ and is always non-positive. Thus, there is no more than one stability switch. 
\end{enumerate}
Hence, due to the continuous dependence on parameters the system has a single stability switch at most. 

\noindent
$\Box$

\medskip

{\bf{Proof of Theorem~\ref{tw5}.}} The quazi-polynomial   has the form
\[
W(\lambda) = \lambda^2 - (a_{11} + a_{22})\lambda\text{e}^{-\lambda \tau} + (a_{11}a_{22} - a_{12}a_{21})\text{e}^{-2\lambda\tau}
 \]
and the characteristic equation $W(\lambda)=0$ is equivalent to
\begin{equation}\label{last}
\lambda^2\text{e}^{2\lambda\tau} - (a_{11} + a_{22})\lambda\text{e}^{\lambda \tau} + a_{11}a_{22} - a_{12}a_{21}=0.
\end{equation}
Let us denote $\lambda \text{e}^{\lambda \tau}=:z$, $a=-(a_{11} + a_{22})=-\text{tr} \mathbf{A}$, $b= a_{11}a_{22} - a_{12}a_{21}=\det \mathbf{A}$. Then Eq.~\eqref{last} reads
\begin{equation}\label{uwik}
z^2+az+b=0,
\end{equation}
and therefore, 
\[
\lambda \text{e}^{\lambda \tau}= \frac{-a-\sqrt{a^2-4ab}}{2} \quad \text{or}
\quad \lambda \text{e}^{\lambda \tau}= \frac{-a+\sqrt{a^2-4ab}}{2} .
\]
Looking for purely imaginary eigenvalues one solves
\[
i\omega_1 \text{e}^{i\omega_1 \tau}= \frac{-a-\sqrt{a^2-4ab}}{2}, \quad
i\omega_2 \text{e}^{i\omega_2 \tau}= \frac{-a+\sqrt{a^2-4ab}}{2}
\]
and obtains two sequences of critical delays at which stability stability switches can occur. 

Now, we need to know the direction of transition of eigenvalues in the complex plane for critical values of delays. We use the theorem of implict function to calculate it. Derivating Eq.~\eqref{uwik} with respect to $\tau$ we obtain
\[
2z\frac{dz}{d\tau}+a\frac{dz}{d\tau}=0 \Longrightarrow \frac{dz}{d\tau}=0 ,
\] 
which yields
\[
\frac{d\lambda}{d\tau}\text{e}^{\lambda \tau}+\lambda\text{e}^{\lambda\tau} 
\left(\frac{d\lambda}{d\tau}\tau+\lambda\right)=0 \Longrightarrow 
\frac{d\lambda}{d\tau}=-\frac{\lambda^2}{1+\lambda \tau} .
\]
Substituting $\lambda=i\omega_k$, $k\in \{1, 2\}$, $\tau=\tau_{cr}$, where $\tau_{cr}$ is any of the critical delays,  one gets 
\[
\frac{d\Re \lambda}{d\tau}\Big|_{\lambda=i\omega_k}=
\frac{\omega^2_k}{1+\omega^2_k\tau_{cr}^2}
\]
which is positive independently on the values of $\omega_k$ and $\tau_{cr}$.
This means that system~\ref{eq:4_del} can have at most one stability switch.

$\Box$

\medskip

\noindent
{\bf{Note to systems with delay in three terms.}}
For Eqs.~(\ref{eq:40}) the quazi-polynomial has the form  
\begin{equation}
W(\lambda) = \lambda^2 - a_{22}\lambda + (-a_{11}\lambda + a_{11}a_{22})\text{e}^{-\lambda\tau} -a_{12}a_{21}\text{e}^{-2\lambda\tau}
\label{eq:655}
\end{equation}
and we can reduce it in three ways:
\begin{enumerate}
\item $a_{11}=0$ yielding $W(\lambda) = \lambda^2 - a_{22}\lambda + -a_{12}a_{21}\text{e}^{-2\lambda\tau}$ for which the auxiliary function $F$ always has one positive zero;
\item $a_{12}a_{21}=0$ yielding $W(\lambda) = \lambda^2 - a_{22}\lambda + (-a_{11}\lambda + a_{11}a_{22})\text{e}^{-\lambda\tau}$ and $F$ has the same property as above;
\item $a_{22}=0$ yielding $W(\lambda) = \lambda^2-a_{11}\lambda\text{e}^{-\lambda\tau} -a_{12}a_{21}\text{e}^{-2\lambda\tau}$ which has the same form as for Eqs.~\eqref{eq:4_del}.
\end{enumerate}
In all three simplifications at most one stability switch is possible.

On the other hand one can rewrite the characteristic equation for $\lambda=i\omega$ in the similar way as when looking for the auxiliary function $F$ and obtain
\begin{equation}
||-\omega^2 - a_{22}i\omega|| =||a_{11}i\omega - a_{11}a_{22} +a_{12}a_{21}(\cos \omega\tau-i\sin\omega\tau)||
\end{equation}
and define the auxiliary function
\[
G_{\tau}(\omega)=\omega^4+ (a_{22}^2-a_{11}^2) \omega^2-a_{11}^2a_{22}^2
-a_{12}^2a_{21}^2+2a_{11}a_{12}a_{21}\sin\omega\tau +2a_{11}a_{22}a_{12}a_{21}\cos\omega\tau
\]
which zeros give the purely imaginary eigenvalues. We see that
\[
G(0)=-(a_{11}a_{22}+a_{12}a_{21})<0 \quad 
\text{and} \quad
\lim_{\omega \to +\infty}G(\omega)=+\infty
\]
implying at least one positive zero of $G$.
We also see that for $\tau=0$ there is $G_0(\omega)=\omega^4+ (a_{22}^2-a_{11}^2) \omega^2-(a_{11}a_{22}+a_{12}a_{21})$ and $G_0$ has one positive zero yielding that for small $\tau>0$ only one positive zero exists. However, looking at the graphs of $\omega^4+ (a_{22}^2-a_{11}^2) \omega^2-a_{11}^2a_{22}^2$ and $a_{12}^2a_{21}^2+2a_{11}a_{12}a_{21}\sin\omega\tau -2a_{11}a_{22}a_{12}a_{21}\cos\omega\tau$ one can easily see that for larger $\tau$ the number of positive zeros increases and there is arbitrary number of positive zeros for sufficiently large $\tau$. This shows that there can be more than one stability switch for some parameters values.

For Eqs.~(\ref{eq:41}) the quazipolynomial has the form 
\begin{equation}
W(\lambda) = \lambda^2 + (-\lambda(a_{11} + a_{22}) - a_{12}a_{21}) \text{e}^{-\lambda\tau} + a_{11}a_{22}\text{e}^{-2\lambda\tau}
\label{eq:655a}
\end{equation}
and we can also reduce it in two ways
\begin{itemize}
\item $a_{11}a_{22}=0$ yielding $W(\lambda) = \lambda^2 + (-\lambda(a_{11} + a_{22}) - a_{12}a_{21}) \text{e}^{-\lambda\tau}$ for which $F$ has one positive zero;
\item  $a_{12}a_{21}=0$ yielding $W(\lambda) = \lambda^2-\lambda(a_{11} + a_{22}) \text{e}^{-\lambda\tau} + a_{11}a_{22}\text{e}^{-2\lambda\tau}$ and $W$ has the same form as for Eqs.~\eqref{eq:4_del}.
\end{itemize}
In both cases only one stability switch is possible. However,  defining the auxiliary function $G$ one can also 
expect multiple stability switches.

\end{document}